\documentclass[a4paper]{pspum-l}
\usepackage{sabbah_augsburg}

\begin{document}

\title[Universal unfoldings of Laurent polynomials and tt$^*$ structures]
{Universal unfoldings of Laurent polynomials\\
and tt$^*$ structures}

\author[C.~Sabbah]{Claude Sabbah}
\address{UMR 7640 du CNRS\\
Centre de Math{\'e}matiques Laurent Schwartz\\
{\'E}cole polytechnique\\
\hbox{F--91128}~Palaiseau cedex\\
France}
\email{sabbah@math.polytechnique.fr}
\urladdr{http://www.math.polytechnique.fr/~sabbah}

\begin{abstract}
This article surveys the relations between harmonic Higgs bundles and Saito structures which lead to tt$^*$ geometry on Frobenius manifolds. We give the main lines of the proof of the existence of a canonical tt$^*$ structure on the base space of the universal unfolding of convenient and nondegenerate Laurent polynomials.
\end{abstract}

\subjclass[2000]{53D45, 32S40, 14C30, 34Mxx}

\keywords{Higgs bundle, Saito structure, tt*-structure, Fourier-Laplace transform, Frobenius manifold, Laurent polynomial}

\maketitle
\tableofcontents

\section*{Introduction}
The notion of a tt$^*$ structure on a holomorphic vector bundle is now understood, after the work of C\ptbl Hertling \cite{Hertling01}, as an enrichment of that of harmonic Higgs bundle previously introduced by N\ptbl Hitchin and C\ptbl Simpson. Given the Higgs field $\Phi$ and the harmonic metric $h$ on the holomorphic bundle $E$ on a complex manifold $M$, the new ingredients needed for a tt$^*$ structure are a real structure on the associated $C^\infty$ bundle $H$, a holomorphic endomorphism $\cU$ of $E$ and a $C^\infty$ endomorphism $\cQ$ of $H$ subject to some compatibility relations. In the following, we relax the condition for $h$ to be positive definite, and only ask that it is Hermitian and nondegenerate. When needed, we will emphasize the positive definite case.

It has been much enlightening in two ways to interpret (\cf \cite{Simpson97}) harmonic Higgs bundles as variations of polarized twistor structures of weight~$0$: firstly, it makes the analogy with variations of Hodge structures more transparent and, secondly, it enables one to do geometry with the external parameter $\hb$ added for this purpose. From this point of view, the relations satisfied by the endomorphisms $\cU$ and $\cQ$ in a tt$^*$ structure appear as expressing the complete integrability of a connection on the twistor bundle.

A nearby (purely holomorphic) notion, that of a Saito structure, has emerged from the work of K\ptbl Saito \cite{KSaito83b} and M\ptbl Saito \cite{MSaito89} as a basic tool to produce Frobenius manifolds from singularities of holomorphic functions. While it is already present in~\cite{Hertling01}, the bridge between these two notions is made more transparent in \S\ref{sec:harmFrob} by the introduction of a potential for the Higgs field.

When the tt$^*$ structure exists on the tangent bundle of $M$, we speak of tt$^*$ geometry, which is a generalization of special geometries on $M$. Of particular interest for us is the case where $M$ is a Frobenius manifold. In such a case, a~Saito structure exists on the tangent bundle together with supplementary symmetry properties, giving rise to a commutative and associative product with unit on $TM$. Adding a~tt$^*$ structure in a compatible way (\ie with the help of a potential for the Higgs field) leads to the structure of \emph{harmonic Frobenius manifold}.

The main result we report here (\cf Theorem \ref{th:main}) is the existence of a canonical harmonic structure on the canonical Frobenius manifold attached to a convenient and nondegenerate Laurent polynomial. Moreover, the corresponding Hermitian form is \emph{positive definite}.

In this survey article, which contains no original result, we first give (\S\S\ref{sec:harmFrob} and~\ref{sec:twstr}) a quick overview of tt$^*$ structures, Saito structures and variations of twistor structures (a more detailed exposition can be found in \cite{Hertling01} and \cite{H-S06}). In \S\ref{sec:Fourier-Laplace}, we explain the Fourier-Laplace method for constructing polarized pure twistor structures starting from a variation of polarized Hodge structure. In \S\ref{sec:laurent}, we show how to apply this technique to the Gauss-Manin connection of a Laurent polynomial, with the help of M\ptbl Saito's mixed Hodge theory \cite{MSaito87}. One can find details for the results of \S\S\ref{sec:Fourier-Laplace} and \ref{sec:laurent} in \cite{Bibi04,Bibi05}, and many other results and applications in \cite{H-S06}.

\subsubsection*{Acknowledgements}
The author thanks Ron Donagi and Katrin Wendland, organizers of the conference ``From tQFT to tt$^*$ and integrability'', for having given him the opportunity to talk about the contents of this article.

\section{Harmonic Frobenius manifolds}\label{sec:harmFrob}
In this section we first recall the notion of a harmonic Higgs bundle, following C\ptbl Simpson \cite{Simpson92}, and we introduce supplementary structures which will be useful for defining harmonic Frobenius manifolds.

\subsection{Harmonic Higgs bundles}\label{subsec:harmbdle}
Let $M$ be a complex manifold and let $E$ be a holomorphic bundle on $M$, equipped with a Hermitian nondegenerate sesquilinear form $h$ (we do not impose at this stage that $h$ is positive definite). We will say that $(E,h)$ is a Hermitian holomorphic bundle. For any operator $P$ acting linearly on~$E$, we will denote by $P^\dag$ its adjoint with respect to $h$.

By a holomorphic \emph{Higgs field}~$\Phi$ on $E$, we mean an $\cO_M$-linear morphism $E\to\Omega^1_M\otimes_{\cO_M}E$ satisfying the ``integrability relation'' $\Phi\wedge\Phi=0$. We then say that $(E,\Phi)$ is a Higgs bundle.

Let $(E,h)$ be a Hermitian holomorphic bundle with Higgs field $\Phi$. Let $H$ be the associated $C^\infty$ bundle, so that $E=\ker d''$, let $D=D'+D''$, with $D''=d''$, be the Chern connection of $h$ and let $\Phi^\dag$ be the $h$-adjoint of $\Phi$. We say that $(E,h,\Phi)$ is a \emph{harmonic Higgs bundle} (or that $h$ is \emph{Hermite-Einstein} with respect to $(E,\Phi)$) if $\VD\defin D+\Phi+\Phi^\dag$ is an integrable connection on~$H$. This is equivalent to a set of relations:
\begin{equation}\label{eq:higgsharmonic}
\left\{
\begin{array}{llll}
(d'')^2&=0,\quad d''(\Phi)&=0,\quad\Phi\wedge\Phi&=0,\\
(D')^2&=0,\quad D'(\Phi^\dag)&=0,\quad\Phi^\dag\wedge\Phi^\dag&=0,\\
D'(\Phi)&=0,\quad d''(\Phi^\dag)&=0,\quad D'd''+{}&d''D'=-(\Phi\wedge\Phi^\dag+\Phi^\dag\wedge\Phi),
\end{array}
\right.
\end{equation}
where the first line is by definition, the second one by $h$-adjunction from the first one, and the third line contains the remaining relations in the integrability condition of~$\VD$. Let us notice that the holomorphic bundle $V=\ker(d''+\Phi^\dag)$ is equipped with a flat holomorphic connection $\Vnablaf$, which is the restriction of $\VD'\defin D'+\Phi$ to $V$.

\subsection{Saito structures on holomorphic bundles}\label{subsec:Saitostr}
The notion of a Saito structure leads, together with a primitive homogeneous section, to the construction of Frobenius manifolds (sometimes one includes the data of the primitive homogeneous section in the Saito structure; we will not do it here). Let us recall it. By a Saito structure (of weight~$0$) on a holomorphic bundle $E$ on $M$, we mean the data of $(E,\nablaf,g,\Phi,\cU,\cV)$ such that $\nablaf$ is a holomorphic connection on $E$, $\Phi$ is a holomorphic Higgs field, $g$ is a nondegenerate symmetric holomorphic bilinear form on $E$, $\cU$ and $\cV$ are endomorphisms of $E$, satisfying (see \eg \cite[\S VI.2.17--2.18]{Bibi00}):
\begin{equation}\label{eq:Saito}
\begin{gathered}
\nablaf^2=0,\quad
\nablaf(\cV)=0,\quad
\Phi\wedge\Phi=0,\quad
[\cU,\Phi]=0\\
\nablaf(\Phi)=0,\quad
\nablaf(\cU)-[\Phi,\cV]+\Phi=0,\\
\nablaf(g)=0,\quad \cV^*+\cV=0\\
\Phi^*=\Phi,\quad \cU^*=\cU.
\end{gathered}
\end{equation}

\subsection{Harmonic Higgs bundles with supplementary structures}\label{subsec:supplstruct}
In this paragraph, we consider supplementary structures on a harmonic Higgs bundle. These supplementary structures produce the data needed to construct a Saito structure (\cf \cite{Hertling01} for more details).

\subsubsection*{Real structure on a Hermitian holomorphic bundle}
Let $(E,h)$ be a Hermitian holomorphic bundle, and let $H$ be the associated $C^\infty$ bundle. By a real structure we mean an \emph{antilinear} isomorphism $\kappa:H\isom H$ such that
\begin{align}
\kappa^2&=\id_H,\label{eq:reel1}\\
h(\kappa\cbbullet,\kappa\cbbullet)&=\ov{h(\cbbullet,\cbbullet)}\quad\text{(equivalently, $h$ is real on $H_\RR\defin\ker(\kappa-\id_H)$)},\label{eq:reel2}\\
D(\kappa)&=0.\label{eq:reel3}
\end{align}

\eqref{eq:reel3} can also be read $D'\kappa=\kappa d''$, or equivalently according to \eqref{eq:reel1}, $\kappa D'=d''\kappa$. Let us consider the nondegenerate complex bilinear form $g(\cbbullet,\cbbullet)\defin h(\cbbullet,\kappa\cbbullet)$. Then, because $h$ is Hermitian, \eqref{eq:reel2} and \eqref{eq:reel3} are respectively equivalent to
\begin{align}
\tag{\ref{eq:reel2}$'$}\label{eq:reel2'}
&\text{$g$ is symmetric},\\
\tag{\ref{eq:reel3}$'$}\label{eq:reel3'}
&\text{$g$ is holomorphic on $E$ and $D'$ is compatible with $g$, \ie $D(g)=0$.}
\end{align}

For any linear operator $P$ on $H$, we denote by $P^*$ the $g$-adjoint of $P$ (which is holomorphic if $P$ is so). Then
\begin{equation}\label{eq:dagstar}
P^*=\kappa P^\dag\kappa.
\end{equation}

\subsubsection*{Potential harmonic Higgs bundles}
Let $(E,h,\Phi)$ be a harmonic Higgs bundle. From \eqref{eq:higgsharmonic}, we have $d''\Phi^\dag=0$. Hence there exists, at least locally, a $C^\infty$ endomorphism $B$ of $E$ such that $d''B=\Phi^\dag$ or, equivalently by setting $A=B^\dag$, $D'A=\Phi$. Let then $(E,h)$ be a Hermitian holomorphic bundle and let $A$ be a $C^\infty$ endomorphism of the associated $C^\infty$ bundle $H$. We say that $(E,h,A)$ is a \emph{potential harmonic Higgs bundle} if
\begin{equation}\label{eq:potentialhiggs}
\text{$(E,h,D'A)$ is a harmonic Higgs bundle.}
\end{equation}
In other words, we assume here the global existence of $A$ and a choice of it. (Except in simple examples as Example \ref{exem:rankone}, it is not expected that $A$ is holomorphic.)

\smallskip
For a potential harmonic Higgs bundle $(E,h,A)$, we set $\Phi=D'A$, which is a Higgs field by definition. We can associate to these data a new connection of type $(1,0)$:
\begin{equation}\label{eq:nablafpotential}
\nablaf\defin D'+[A^\dag,\Phi].
\end{equation}
From \eqref{eq:higgsharmonic} we get that \emph{$\nablaf$ is holomorphic}, \ie $d''\nablaf+\nablaf d''=0$.

\subsubsection*{Real (potential) harmonic Higgs bundles}
Let us consider $(E,h,\Phi,\kappa)$ (resp.\ $(E,h,A,\kappa)$), where $(E,h,\kappa)$ is a real Hermitian holomorphic bundle and $(E,h,\Phi)$ (\resp $(E,h,A)$) is a (\resp potential) harmonic Higgs bundle. We say that $(E,h,\Phi,\kappa)$ is a real harmonic Higgs bundle if the following compatibility relation is satisfied:
\begin{equation}\label{eq:realharmonicHiggs}
\Phi^*=\Phi\quad\text{(\resp $A^*=A$)}.
\end{equation}
According to \eqref{eq:dagstar}, this is equivalent to $\Phi^\dag=\kappa\Phi\kappa$ (\resp $A^\dag=\kappa A\kappa$) or, equivalently, to $\VD(\kappa)=0$, \ie $\kappa$ is a real structure on the local system $\ker\Vnablaf$. In other words, if we denote, for any $\hbo\in\CC^*$, by $(V_{(\hbo)},\Vhbonablaf)$ the flat bundle $\big(\ker (D''+\nobreak\hbo\Phi^\dag),D'+\nobreak\Phi/\hbo\big)$, then the conditions \eqref{eq:reel1}, \eqref{eq:reel3} and \eqref{eq:realharmonicHiggs} on $\kappa$ can be rephrased by saying that $\kappa$ is a real structure on each of the locally constant sheaves $\ker\Vhbonablaf$, $\hbo=\pm1$ (or equivalently, $|\hbo|=1$).

Let us notice that, when $\Phi=D'A$, the condition $A^*=A$ implies $\Phi^*=\Phi$, as we have $D'(A^*)=(D'A)^*$ (because $D'(g)=0$). Conversely, If $\Phi=D'A$ and $\Phi=\Phi^*$, we also have $\Phi=D'(A+A^*)/2$, so we can assume that $A=A^*$. Moreover, with such an assumption, $[A^\dag,\Phi]^*+[A^\dag,\Phi]=0$, so, according to \eqref{eq:reel3'},
\begin{equation}\label{eq:nablafg}
\nablaf(g)=0.
\end{equation}
Lastly, from $D(\Phi)=0$ given by \eqref{eq:higgsharmonic} we also deduce
\begin{equation}\label{eq:nablafPhi}
\nablaf(\Phi)=0.
\end{equation}

\subsubsection*{Integrable (potential) harmonic Higgs bundles}
Let $(E,h,\Phi)$ be a harmonic Higgs bundle. We say that it is \emph{integrable}\footnote{This terminology will be justified in \S\ref{subsec:integrability}.} if there exist two endomorphisms $\cU$ and $\cQ$ of $H$ such that
\begin{align}
\text{$\cU$ is holomorphic, \ie }d''(\cU)&=0,\label{eq:higgsint1}\\
\cQ^\dag&=\cQ,\label{eq:higgsint2}\\
[\Phi,\cU]&=0,\label{eq:higgsint3}\\
D'(\cU)-[\Phi,\cQ]+\Phi&=0,\label{eq:higgsint4}\\
D'(\cQ)+[\Phi,\cU^\dag]&=0\label{eq:higgsint5}.
\end{align}
Let us note that $(\cU,\cQ+\lambda\id)$ satisfy the same equations for any $\lambda\in\RR$.

If $(E,h,\Phi)$ is moreover potential (\ie $\Phi=D'A$), then, in analogy with \eqref{eq:nablafpotential}, we~set
\begin{equation}\label{eq:cVpotential}
\cV=\cQ-[A^\dag,\cU].
\end{equation}
The equation $\eqref{eq:higgsint5}^\dag$ obtained from \eqref{eq:higgsint5} by adjunction can be written $d''(\cV)=0$, showing that $\cV$ is then a holomorphic endomorphism of $E$.

Let us also notice that \eqref{eq:higgsint4} can as well be written as
\begin{equation}\label{eq:nablafcU}
\nablaf(\cU)-[\Phi,\cV]+\Phi=0.
\end{equation}

\subsubsection*{tt$^*$ bundles}
We now end with all the structures put together, that is, we consider t-uple $(E,h,\Phi,\kappa,\cU,\cQ)$ such that $(E,h,\kappa)$ is a real Hermitian holomorphic bundle (so \eqref{eq:reel1}--\eqref{eq:reel3} are satisfied), $(E,h,\Phi,\cU,\cQ)$ is an integrable harmonic Higgs bundle (\ie \eqref{eq:potentialhiggs} and \eqref{eq:higgsint1}--\eqref{eq:higgsint5} are satisfied), and moreover
\begin{align}
\cU^*&=\cU,\label{eq:realintHiggscU}\\
\cQ^*+\cQ&=0.\label{eq:realintHiggscQ}
\end{align}
Let us note that \eqref{eq:higgsint2} and \eqref{eq:realintHiggscQ} imply
\begin{equation}\label{eq:cQkappa}
\text{$\cQ$ is purely imaginary, \ie }\cQ\kappa=-\kappa\cQ.
\end{equation}

One defines potential tt$^*$ bundles similarly. For such an object,
\begin{equation}\label{eq:realintHiggscV}
\cV^*+\cV=0.
\end{equation}

\subsubsection*{Data for a Saito structure}
Given potential tt$^*$ bundle $(E,h,A,\kappa,\cU,\cQ)$, we have associated to it a t-uple $(E,\nablaf,g,\Phi,\cU,\cV)$ satisfying the relations \eqref{eq:Saito} except possibly
\begin{equation}\label{eq:harmpreSaito}
\nablaf^2=0,\quad\nablaf(\cV)=0.
\end{equation}
(Indeed, these relations were derived in \eqref{eq:reel3'}, \eqref{eq:realharmonicHiggs}, \eqref{eq:nablafg}, \eqref{eq:nablafPhi}, \eqref{eq:nablafcU}, \eqref{eq:realintHiggscU} and \eqref{eq:realintHiggscV}.)

\subsection{Saito structures with supplementary structures}\label{subsec:Saitosuppl}
Starting now from a Saito structure $(E,\nablaf,g,\Phi,\cU,\cV)$ on $E$ as in \S\ref{subsec:Saitostr}, we will add supplementary structures $\kappa$ and $A$ and try to reconstruct a tt$^*$ bundle.

\subsubsection*{Real structure}
We still denote by $H$ the $C^\infty$ bundle associated to $E$. We say that an antilinear isomorphism $\kappa:H\to H$ is a real structure if $\kappa^2=\id_H$. We define the nondegenerate sesquilinear form $h$ by $h(\cbbullet,\kappa\cbbullet)=g(\cbbullet,\cbbullet)$. We impose the following compatibility relations:
\begin{itemize}
\item
$g$ is real on $H_\RR\defin\ker(\kappa-\id_H)$,
\item
$D'(g)=0$, where $D=D'+d''$ is the Chern connection of $h$.
\end{itemize}
We then speak of a \emph{real Saito structure}. We say that the real Saito structure $(E,\nablaf,g,\Phi,\cU,\cV,\kappa)$ is \emph{harmonic} if $(E,h,\Phi)$ is a harmonic Higgs bundle.

\subsubsection*{Potential for a Saito structure}
Let $(E,\nablaf,g,\Phi,\cU,\cV)$ be a Saito structure and let $A^\dag$ be a $C^\infty$ endomorphism of $E$ such that $A^{\dag*}=A^\dag$ (here, $^\dag$ makes no reference to any adjunction). We then set $\Phi^\dag=d''A^\dag$, which satisfies thus $\Phi^{\dag*}=\Phi^\dag$. Let us then define
\begin{align*}
D'&=\nablaf-[A^\dag,\Phi],\\
\cQ&=\cV+[A^\dag,\cU].
\end{align*}
We notice that, with this definition of $D'$ and $\Phi^\dag$, the third line of \eqref{eq:higgsharmonic} is satisfied and that $D'(g)=0$ (hence $D(g)=0$, if we set $D=D'+d''$). We also have $\cQ^*+\cQ=0$, \eqref{eq:higgsint4} is fulfilled as well as \eqref{eq:higgsint5}$^\dag$.

We say that $A^\dag$ is a harmonic potential for the Saito structure $(E,\nablaf,g,\Phi,\cU,\cV)$ if the second line of \eqref{eq:higgsharmonic} is satisfied.

\subsubsection*{Harmonic potential real Saito structure}
Let $(E,\nablaf,g,\Phi,\cU,\cV)$ be a Saito structure equipped with a harmonic real structure $\kappa$ (and associated harmonic Hermitian form $h$) and a harmonic potential~$A^\dag$. We say that the two structures are compatible if the following properties are fulfilled:
\begin{enumerate}
\item
$h$ and $A^\dag$ define the same connection $D'$,
\item
the $h$-adjoint of $\Phi$ is $\Phi^\dag\defin d''A^\dag$,
\item
$\cQ^\dag=\cQ$.
\end{enumerate}
We then speak of a harmonic potential real Saito structure $(E,\nablaf,g,\Phi,\cU,\cV,\kappa,A^\dag)$.

\begin{proposition}\label{prop:SaitoHiggs}
There is a one-to-one correspondence between harmonic potential real Saito structures $(E,\nablaf,g,\Phi,\cU,\cV,\kappa,A^\dag)$ and potential tt$^*$ bundles $(E,h,A,\kappa,\cU,\cQ)$ satisfying \eqref{eq:harmpreSaito}.\qed
\end{proposition}

\subsection{Harmonic Frobenius manifolds}
A structure of Frobenius manifold (or Saito manifold) on a complex manifold $M$ consists in giving a Saito structure on the tangent bundle $M$ with the supplementary conditions that
\begin{align}
&\text{$\nablaf$ is torsion free},\label{eq:frob1}\\
&\text{$\Phi$ is symmetric},\label{eq:frob2}\\
&\exists e\in\Gamma(M,TM),\quad\text{$\nablaf(e)=0$, $e$ is an eigenvector of $\cV$ and $\Phi_e=-\id_{TM}$.}\label{eq:frob3}
\end{align}

Let us recall that a Higgs field on $TM$ defines a product $\star$ on $TM$ by setting $\xi\star\eta=-\Phi_\xi\eta$ and that saying $\Phi$ is symmetric means that $\star$ is commutative (and also associative because of the Higgs condition). The vector field $e$ is then a unit for~$\star$. The \emph{Euler field} $\gE$ is defined as $\cU(e)$, and $\cU$ is identified with the multiplication by~$\gE$. If $c\in\CC$ is defined by $\cV e=ce$, we then have $\nablaf\gE=\cV+(1-c)\id$, so that $\nablaf_e\gE=e$. Moreover, we have $\Lie_\gE(g)=2(1-c)$ (\cf \eg \cite[\S VII.1.21]{Bibi00}).

\begin{definition}[Harmonic Frobenius manifold]
A structure of \emph{harmonic Frobenius manifold} (or \emph{harmonic Saito manifold}) on a complex manifold $M$ consists in giving a \emph{harmonic potential real Saito structure} on the tangent bundle $M$ with the supplementary conditions \eqref{eq:frob1}--\eqref{eq:frob3} and such that $c\in\RR$.
\end{definition}

The new data added to a Frobenius manifold in order to make it harmonic consist of a real structure $\kappa$ on $TM$ and a potential $A^\dag$ (or $A$), satisfying the relations described in \S\ref{subsec:Saitosuppl}.

\begin{proposition}
A harmonic Frobenius manifold is a manifold with a CDV-structure, in the sense of \cite[Def\ptbl1.2]{Hertling01}.
\end{proposition}

\begin{proof}
Firstly, (1.5) in \loccit means $\kappa\Phi\kappa=\Phi^\dag$, which is equivalent to $\Phi^*=\Phi$ (and is true by assumption, \cf \eqref{eq:Saito}).

As $\nablaf e=0$ and $\nablaf$ is torsion free, we have $\nablaf_e=\Lie_e$. As $\Phi_e=-\id$, \eqref{eq:nablafpotential} implies then $\Lie_e=\nablaf_e=D'_e$. Therefore $\Lie_e(h)=0$, so \cite[(1.6)]{Hertling01} is satisfied.

We will now show that, if we denote by $c$ the eigenvalue of $\cV$ with respect to~$e$ as above, we have
\begin{equation}\label{eq:cQDE}
\cQ=D_\gE-\Lie_\gE-(1-c)\id.
\end{equation}
Let us first remark that, as $\gE$ is holomorphic, $d''_\gE$ acting on $TM$ is zero, so that we can replace $D_\gE$ with $D'_\gE$ in \eqref{eq:cQDE}. Then, by \eqref{eq:cVpotential}, $\cQ=\cV+[A^\dag,\cU]=\nablaf\gE+[A^\dag,\cU]-(1-c)\id$. But on the other hand, using the torsion freeness of $\nablaf$,
\[
\nablaf\gE=\nablaf_\gE-\Lie_\gE=D'_\gE+[A^\dag,\Phi_\gE]-\Lie_\gE=D'_\gE-[A^\dag,\cU]-\Lie_\gE,
\]
hence \eqref{eq:cQDE}. In other words, \cite[(1.8)]{Hertling01} is satisfied.

We use \cite[Rem\ptbl3.6]{Takahashi04} to remark that $\cQ^\dag=\cQ$ and $c\in\RR$ imply $\Lie_{\gE-\ov\gE}(h)=0$, so \cite[(1.7)]{Hertling01} is satisfied.

Lastly, the integrability property \cite[(1.9)]{Hertling01} is a consequence of \eqref{eq:higgsint1}--\eqref{eq:higgsint5} and the adjoint relations.
\end{proof}

Following the method of K\ptbl Saito, the infinitesimal period mapping with respect to a primitive section (\cf \eg \cite[Th\ptbl5.12]{Hertling01}, \cite[Chap\ptbl VII]{Bibi01c}) gives (\cf \cite[Th\ptbl5.15]{Hertling01}):

\begin{corollaire}\label{cor:harmonic}
Let $(E,\nablaf,g,\Phi,\cU,\cV,\kappa,A^\dag)$ be a harmonic potential real Saito structure on $E$ and let $\omega$ be a $\nablaf$-horizontal section of $E$ which is primitive and real homogeneous, \ie such that
\begin{enumerate}
\item
$\varphi_\omega(\cbbullet)\defin-\Phi_\bbullet(\omega):TM\to E$ is an isomorphism,
\item
$\omega$ is an eigenvector of $\cV$ with real eigenvalue.
\end{enumerate}
Then $\varphi_\omega$ equips $M$ with the structure of a harmonic Frobenius manifold.\qed
\end{corollaire}

\subsection{Examples}

\begin{exemple}[Compact K\"ahler manifolds, after K\ptbl Corlette and C\ptbl Simpson]\label{exem:corlette-simpson}
The main source of examples, when $M$ is compact K\"ahler, comes from the fundamental result of K\ptbl Corlette \cite{Corlette88} and C\ptbl Simpson \cite{Simpson88} saying that a flat bundle $(V,\Vnablaf)$ admits a harmonic \emph{metric} (\ie positive definite) if and only if the local system $\ker\Vnablaf$ is \emph{semi-simple}. The metric is then uniquely determined (up to a positive constant) on any simple summand. Similarly, a Hermite-Einstein metric $h$ exists for $(E,\Phi)$ iff $(E,\Phi)$ is \emph{poly-stable with vanishing Chern classes}.

In general, such a metric is not determined explicitly, but is given by a nontrivial global existence result.
\end{exemple}

\begin{exemple}[Quasi-projective varieties, after T\ptbl Mochizuki]\label{exem:mochizuki}
The results of Corlette and Simpson have been generalized by C\ptbl Simpson in dimension one \cite{Simpson90} and then by T\ptbl Mochizuki \cite{Mochizuki07} in arbitrary dimension. If $M$ is a smooth quasi-projective variety and $\wt M$ is a smooth projective variety containing $M$ as a Zariski dense subset and such that $D\defin\wt M\moins M$ is a normally crossing divisor in $\wt M$, then one can introduce (\cf \cite{Simpson90} if $\dim M=1$) the notion of tameness for a Higgs field on a holomorphic vector bundle $E$ on $M$: the eigenvalues of the Higgs field should have at most logarithmic growth along $D$. Similarly, given a flat bundle $(V,\Vnablaf)$ with metric $h$ on $M$, tameness means that the $h$-norm of horizontal sections has at most a moderate growth along $D$. One of the main results of \cite{Mochizuki07} is a one-to-one correspondence between semi-simple local systems on $M$ and harmonic flat bundles $(V,\Vnablaf,h)$ for which the metric is tame and the eigenvalues of the residue of the Higgs field along $D$ are purely imaginary. Here also, the metric is not determined explicitly, but is given by a nontrivial global existence result.
\end{exemple}

\begin{exemple}[Variations of complex Hodge structures of weight $0$]\label{exem:VHS}
Let $H$ be a $C^\infty$ vector bundle on $M$, equipped with a flat connection $\VD$ and a decomposition $H=\oplus_{p\in\ZZ}H^p$ by $C^\infty$ subbundles. We assume that Griffiths transversality relations hold:
\[
\VD' H^p\subset (H^p\oplus H^{p-1})\otimes_{\cO_M}\Omega^1_M,\quad\VD''H^p\subset (H^p\oplus H^{p+1})\otimes_{\cO_{\ov M}}\Omega^1_{\ov M}.
\]
We denote by $D_{|H^p}$ the composition of $\VD_{|H^p}$ with the projection to $H^p$, by $\Phi_{|H^p}$ the composition of $\VD'_{|H^p}$ with the projection to $H^{p-1}$ and by $\Phi^\dag$ that of $\VD''_{|H^p}$ with the projection to $H^{p+1}$. We then set $D=\oplus_pD_{|H^p}$, $\Phi=\oplus_p\Phi_{|H^p}$ and $\Phi^\dag=\oplus_p\Phi^\dag_{|H^p}$. Considering the holomorphic type together with the degree with respect to the decomposition, we find that the relations \eqref{eq:higgsharmonic} are satisfied because of the flatness of $\VD$.

Assume that we are given a nondegenerate Hermitian form $k$ such that \hbox{$\VD(k)=0$} and the decomposition $H=\oplus_{p\in\ZZ}H^p$ is $k$-orthogonal. Consider the nondegenerate Hermitian form $h=\oplus_p(-1)^pk_{|H^p}$. Then $D(h)=0$ and $\Phi^\dag$ is the $h$-adjoint of~$\Phi$. In such a case, we say that $(H=\oplus_pH^p,\VD,k)$ is a variation of complex Hodge structure of weight $0$. In particular, $(H,D'',h,\Phi)$ is a harmonic Higgs bundle. Let us show that it is integrable. We set $\cQ=\oplus_pp\id_{H^p}$ and $\cU=0$. By definition, $D$~is compatible with the decomposition, hence $D(\cQ)=0$ and, as $p$ is real, we have $\cQ^\dag=\cQ$. Lastly, we have $[\Phi,\cQ]=\Phi$.

By a real structure $\kappa$, we mean an anti-linear involution $\kappa:H\to H$ which is $\VD$-horizontal (defining thus a real structure on the local system $\ker\VD$) such that $\kappa(H^p)=H^{-p}$ for any $p$. Then $D(\kappa)=0$ and $\Phi^\dag=\kappa\Phi\kappa$. The previous data define thus a tt$^*$ bundle.
\end{exemple}

\begin{exemple}[The rank-one case on a punctured disc]\label{exem:rankone}
In order to have a concrete example at hand, we will give details on the rank-one case on the punctured unit disc $\Delta^*$ (having coordinate $t$), with a positive definite Hermitian form $h$ (\ie a Hermitian metric). Although very simple, the rank-one case gives indication on the relations between the eigenvalues of the Higgs field and that of the monodromy of the integrable connection $\Vnablaf$. Let us first introduce some notation.

In the following, we denote by $\gB$ the set of equivalence classes of affine forms $\gb:\hb\mto u\hb+v$ modulo the addition of a purely imaginary complex number: we have $\gb\equiv\gb'\ssi (u=u') \text{ and } (v-v'\in i\RR)$. Given an affine form $\gb$, we consider the two functions\footnote{In \cite{Mochizuki07}, Mochizuki uses the notation $-\mathfrak{p},\mathfrak{e}$.}
\begin{align*}
\hb&\mto b_\hb\defin\reel\gb(\hb),\\
\hb&\mto\beta_\hb\defin\hb\Retw\gb(\hb),
\end{align*}
where, if $\gb(\hb)=u\hb+v$, we set\footnote{Anticipating the notation introduced in \S\ref{sec:twstr}, let us notice that $\Retw\gb(\hb)=\frac12(\gb(\hb)+\ovv{\gb(\hb)})$.}
\[
\hb\Retw\gb(\hb)=\tfrac12\hb\big[u\hb+v-\ov u\hbm+\ov v\big]=\tfrac12 u\hb^2+(\reel v)\hb-\tfrac12\ov u.
\]
Therefore, the two functions $b_\bbullet$ and $\beta_\bbullet$ only depend on the class of $\gb$ in $\gB$. Moreover, the function $\beta_\bbullet$ uniquely determines the class of $\gb\in\gB$: indeed, the coefficients of~$\beta_\bbullet$ give $u$ and $\reel v$.

\refstepcounter{equation}\label{stat:rankone}
\par\smallskip
\noindent\eqref{stat:rankone}\enspace
\emph{The set of isomorphism classes of rank-one harmonic Higgs bundles (or harmonic flat bundles) on $\Delta^*$ is in one-to-one correspondence with the set of pairs $(\psi,\gb)$, with $\psi\in\cO(\Delta^*)$ and $\psi dt/t$ having no residue and $\gb\in\gB\bmod\ZZ$. Moreover, meromorphic Higgs bundles correspond to pairs $(\psi,\gb)$ with $\psi$ meromorphic at~$0$.}

\begin{proof}
We will do the proof for harmonic Higgs bundles. The case of harmonic flat bundles is similar and the proof will show the 1-1 correspondence between both. Let us start with a rank-one harmonic Higgs bundle $(H,D'',\Phi,h)$. Let $e'$ be some holomorphic frame of $E=\ker D''$. The Higgs field $\Phi$ can be written as
\[
\Phi e'=\varphi(t)\,e'\,dt,\quad \text{$\varphi(t)$ holomorphic on $\Delta^*$.}
\]
Let us set $\varphi(t)=\partial_t\psi(t)+\beta_0/t$ with $\psi\in\cO(\Delta^*)$ and $\psi dt/t$ having no residue at $0$ and $\beta_0\in\CC$. Let us note here that $\varphi$ is meromorphic iff $\psi$ is so. Let us also set $\norme{e'}_h=\exp(\lambda(t))$ for some $C^\infty$ real function $\lambda$ on $\Delta^*$. Then we have $D'e'=2\partial_t\lambda\cdot e'\, dt$ and the harmonicity condition on the Higgs bundle is that $\lambda$ is \emph{harmonic} on $\Delta^*$. Indeed, in the orthonormal frame $\epsilon=\exp(-\lambda(t))\cdot e'$, we have
\begin{equation}\label{eq:DEeps}
D'\epsilon=\partial_t\lambda\,\epsilon\,dt,\quad D''\epsilon=-\partial_{\ov t}\lambda\,\epsilon\,d\ov t,\quad \Phi\epsilon= \varphi(t)\,\epsilon\,dt,\quad \Phi^\dag\epsilon=\ov \varphi(t)\,\epsilon\,d\ov t.
\end{equation}
Therefore,
\begin{equation}\label{eq:DVeps}
\VD'\epsilon=(\varphi(t)+\partial_t\lambda)\,\epsilon\,dt,\quad \VD''\epsilon=(\ov \varphi(t)-\partial_{\ov t}\lambda)\,\epsilon\,d\ov t.
\end{equation}
The flatness of $\VD$ is then equivalent to $2\partial_t\ov\partial_t\lambda=0$.

If $\lambda$ is harmonic on $\Delta^*$, there exists (by working on the universal cover of $\Delta^*$) a holomorphic function $\mu$ on $\Delta^*$ and a real number $c$, such that $\lambda(t)=2(\reel\mu(t)+\nobreak c\log\mt)$. We have $\partial_t\lambda=\partial_t\mu+c/t$ and $\partial_{\ov t}\lambda=\ov{\partial_t\lambda}$. Let us then replace $e'$ with $e=\exp(-2\mu(t))\cdot e'$, which has $h$-norm $\mt^{2c}$. So we can assume from the beginning that $\lambda=b_0\log\mt$, $b_0\in\RR$. Therefore, $\epsilon=\mt^{-b_0}e$ and we have
\[
D'e=\frac{b_0}{t}\,e\,dt,\quad \Phi e=\Big(\partial_t\psi(t)+\frac{\beta_0}{t}\Big)e\,dt.
\]

Now, \eqref{eq:DEeps} and \eqref{eq:DVeps} can be rewritten as
\begin{gather}
\begin{aligned}
D'\epsilon&=\frac{b_0}{2t}\,\epsilon\,dt,& D''\epsilon&=-\frac{b_0}{2\ov t}\,\epsilon\,d\ov t,\\
\Phi\epsilon&= \Big(\partial_t\psi(t)+\frac{\beta_0}{t}\Big)\epsilon\,dt,&\Phi^\dag\epsilon&=\Big(\ov{\partial_t\psi(t)+\frac{\beta_0}{t}}\Big)\epsilon\,d\ov t
\end{aligned}
\\
\VD'\epsilon=\Big(\partial_t\psi(t)+\frac{b_0+2\beta_0}{2t}\Big)\epsilon\,dt,\quad \VD''\epsilon=\Big(\ov{\partial_t\psi(t)+\frac{2\beta_0-b_0}{2t}}\Big)\epsilon\,d\ov t.
\end{gather}

We will consider the affine form
\[
\gb(\hb)=-2\ov\beta_0\hb+b_0.
\]
We will use the notation $b_1=(\reel\gb)(1)=b_0-2\reel \beta_0$, $\beta_1=(\hb\Retw\gb)(1)=b_0+\beta_0-\ov \beta_0$. Let us consider the frame
\[
v=e^{\psi-\ov \psi}\mt^{\gb(1)}\epsilon.
\]
It is $\VD''$-holomorphic and has $h$-norm $\mt^{b_1}$. We have
\[
\Vnablaf v=\Big(2\partial_t\psi+\frac{\beta_1}{t}\Big)v\,dt.
\]

Let us end the proof of \eqref{stat:rankone}. Let $e'$ be another holomorphic frame of~$E$ such that $\mt^{-b'_0}e'$ has norm one. Then $e'=\nu(t)e$ for some holomorphic function $\nu(t)$ on $\Delta^*$ having moderate growth, which is thus meromorphic. We hence have $b'_0-b_0=k\in\ZZ$ and $e'=\nu t^k e$ with $\module{\nu}=1$. The Higgs field has the same expression in the bases $e$ and $e'$, fixing thus $\beta_0$ and $\psi$.
\end{proof}

\begin{remarque}
The reason for defining $\gB$ as equivalence classes modulo $i\RR$ is to be able to conclude that both categories of rank-one harmonic Higgs bundles and rank-one harmonic flat bundles are equivalent, even if the representative $\gb=u\hb+v$ chosen in each case is not the same. In the Higgs case, we choose $\im v=0$, while in the flat case we choose $\im v=-\im u$.
\end{remarque}

In this example, tameness (as mentioned in Example \ref{exem:mochizuki}) means $\psi$ holomorphic on $\Delta$, while the condition of being purely imaginary means $\beta_0\in i\RR$, so $b_1=b_0$.

A real structure exists if and only if the residue of $\Vhbonablaf$ with $\hbo=\pm1$ belongs to $\hZZ$. This reduces to $\beta_1$ and $-\beta_{-1}\in\hZZ$, that is, $b_0\in\hZZ$ and $\beta_0\in\RR$.

A potential $A$ is now a $C^\infty$ function such that $\partial_tA=\partial_t\psi+\beta_0/t$. One must have $A=\psi+\beta_0\log|t|^2+\text{anti-holomorphic}$. However, as $A^\dag$ only appears through a commutator, it is not needed in rank one, and one has $\nablaf=D'$.

If $\beta_0\neq0$, integrability only holds on simply connected open sets of the punctured disc: $\cQ$ must be a real constant and $\cU$ has to be a holomorphic function which satisfies $\cU'(t)=-(\partial_t\psi+\beta_0/t)$. In order to have integrability on the whole punctured disc while keeping a singularity at the origin, it is necessary to accept a nontame metric, \ie $\psi$ meromorphic and nonholomorphic.
\end{exemple}

\section{Integrable variations of twistor structures}\label{sec:twstr}
In this section, we explain a criterion, due to C\ptbl Hertling, for proving the harmonicity of a Frobenius manifold. We first recall the introduction of an extra parameter to explain either the relations defining a harmonic Higgs bundle or that defining a Saito structure.

All along this section, $M$ will denote a complex manifold, $\ov M$ the conjugate manifold (with structure sheaf $\ov{\cO_M}$), and $M_\RR$ will denote the underlying real-analytic or $C^\infty$-manifold. We will denote by $\PP^1$ the Riemann sphere, covered by the two affine charts $\simeq\Afu$ with coordinate $\hb$ and $1/\hb$, and by $p:M\times\PP^1\to M$ the projection.

The coordinate $\hb$ being fixed, we denote by $\bS$ the circle $|\hb|=1$, by $\Omega_0$ an open neighbourhood of the closed disc $\Delta_0\defin\{|\hb|\leq1\}$ and by $\Omega_\infty$ an open neighbourhood of the closed disc $\Delta_\infty\defin\{|\hb|\geq1\}$.

We will denote by $\sigma:\PP^1\to\ov{\PP^1}$ the anti-holomorphic involution $\hb\mto-1/\ov\hb$. We assume that $\Omega_\infty=\sigma(\Omega_0)$. We denote by $\iota:\PP^1\to\PP^1$ the holomorphic involution $\hb\mto-\hb$.

\subsection{Families of vector bundles on the Riemann sphere}

Let $\cH''$ be a holomorphic vector bundle on $M\times\Omega_0$. Then $\ov{\cH''}$ is a holomorphic bundle on the conjugate manifold $\ov M\times\ov{\Omega_0}$ and $\sigma^*\ov{\cH''}$ is a holomorphic bundle on $\ov M\times\Omega_\infty$ (\ie is an anti-holomorphic family of holomorphic bundles on $\Omega_\infty$). We will set $\ovv{\cH''}\defin\sigma^*\ov{\cH''}$.

By a \emph{real-analytic} (\resp $C^\infty$) \emph{family} of holomorphic vector bundles on $\PP^1$ parametrized by $M_\RR$ we will mean the data of a triple $(\cH',\cH'',\cCS)$ consisting of holomorphic vector bundle $\cH',\cH''$ on $M\times\Omega_0$ and a nondegenerate $\cO_{M\times\bS}\otimes_{\cO_\bS}\cO_{\ov M\times\bS}$-linear morphism
\[
\cCS:\cHS'\otimes_{\cO_\bS}\ovv{\cHS''}\to\cC^{*,\an}_{M_\RR\times\bS},
\]
where $\cC^{*,\an}_{M_\RR\times\bS}$ is either the sheaf of real-analytic functions on $M_\RR\times\bS$ if $*=\Ran$, or the sheaf of $C^\infty$ functions on $M_\RR\times\bS$ which are real analytic with respect to $\hb\in\bS$ if $*=\infty$. The nondegeneracy condition means that $\cCS$ defines $C^{*,\an}$-gluing between the dual $\cH^{\prime\vee}$ of $\cH'$ and $\ovv{\cH''}$, giving rise to a $\cC^{*,\an}_{M_\RR\times\PP^1}$-locally free sheaf of finite rank that we denote by $\cH$, where $\cC^{*,\an}_{M_\RR\times\PP^1}$ denotes the sheaf of $C^\infty$ or real-analytic functions on $M_\RR\times\PP^1$ \emph{which are holomorphic with respect to~$\hb$}.

\begin{remarque}
We could have started from $\cC^{*,\an}_{M_\RR\times\Omega_0}$-locally free sheaves (\cf \eqref{eq:cC} below), but we insist here that these sheaves come from holomorphic bundles, even if the result $\cH$ of the gluing is only a $\cC^{*,\an}_{M_\RR\times\PP^1}$-locally free sheaf.

On the other hand, using a pairing instead of a gluing is not essential here, as we can replace $\cH'$ with its dual bundle, and is mainly motivated by \eqref{eq:hS}, as well as a better behaviour near singularities, when one develops the theory with singularities (twistor $\cD$-modules, \cf \cite{Bibi01c}).
\end{remarque}

In order to give a uniform result on holomorphic, real-analytic or $C^\infty$ families of holomorphic bundles on $\PP^1$, we will give a similar presentation for holomorphic families. Let us denote by $\wt\sigma:\PP^1\to\PP^1$ the involution $\hb\mto-1/\hb$. Given a holomorphic bundle $\cH''$ on $M\times\Omega_0$, we now denote by $\ovv{\cH''}$ the holomorphic bundle $\wt\sigma^*\cH''$.

By a \emph{holomorphic family} of vector bundles on $\PP^1$ parametrized by $M$, we will mean the data of a triple $(\cH',\cH'',\cCS)$, where $\cH'$ (\resp $\cH''$) is a holomorphic vector bundle on $M\times\Omega_0$ and $\cCS:\cHS'\otimes_{\cO_{M\times\bS}}\ovv{\cHS''}\to\cO_{M\times\bS}$ is nondegenerate $\cO_{M\times\bS}$-bilinear pairing between the sheaf-theoretic restrictions to $M\times\bS$, defining a gluing between $\cH^{\prime\vee}$ and $\ovv{\cH''}$, hence a holomorphic vector bundle $\cH$ on $M\times\PP^1$.

Let us notice that a real analytic family as above, parametrized by the real analytic variety $M_\RR$, can be extended to the germ of complexification $(M_\RR)_\CC$ of $M_\RR$ as a holomorphic family.

\begin{proposition}\label{prop:familyP1}
Let $(\cH',\cH'',\cCS)$ be a holomorphic (\resp $C^{\Ran,\an}$) family of holomorphic vector bundles on $\PP^1$ parametrized by $M$ (\resp $M_\RR$).
\begin{enumerate}
\item\label{prop:familyP11}
If $\cH_{|\{x\}\times\PP^1}$ is the trivial bundle on $\PP^1$ for any $x\in M$, then $H\defin p_*\cH$ is a holomorphic (\resp $C^{\Ran}$) vector bundle on $M$ (\resp $M_\RR$) and the natural morphism $p^*H=p^*p_*\cH\to\cH$ is an isomorphism.
\item\label{prop:familyP12}
If $M$ is connected, the subset $\Theta$ of $M$ defined by the property that $\cH_{|\{x\}\times\PP^1}$ is the trivial bundle on $\PP^1$ iff $x\in M\moins\Theta$ is either empty, or equal to $M$, or locally defined by a single holomorphic (\resp real-analytic) equation.
\end{enumerate}
Moreover, the first point also holds for $C^\infty$ families.
\end{proposition}

\begin{proof}
In the case of holomorphic families, one can use an argument on direct images of coherent sheaves, characterizing $\Theta$ as the support of the coherent sheaf $R^1p_*\cH(-1)$ (\cf \eg \cite[Th\ptbl I.5.3]{Bibi00}). The case of real analytic families is obtained by applying the case of holomorphic families to the complexified family.

For the first point in the case of $C^\infty$ families, one can adapt the proof of \cite[\S4]{Malgrange83db}, which also applies to the holomorphic or real analytic case. One first notices that this is a local result on $M$, so that we can assume that $M$ is a polydisc in $\CC^n$. Then, applying Grauert's theorem (or a simple variant of it, as \cite[Cor\ptbl2.17]{Leiterer90}), we can assume that $\cH'$ and $\cH''$ are trivial holomorphic vector bundles. Let us denote by $d$ the rank of $\cH',\cH''$. Fixing the trivializations, the pairing $\cCS$ is now regarded as an invertible $d\times d$-matrix $C_x(\hb)$ defined on $M\times\{1/r<\hb<r\}$ with $r>1$ (up to shrinking $M$ once more), which is holomorphic with respect to $\hb$ and $C^\infty$ with respect to $x$.

The bundle $\cH_{|\{x\}\times\PP^1}$ is trivial if, up to reducing $r>1$, there exist holomorphic invertible matrices $\hb\mto\Sigma'_x(\hb),\Sigma''_x(\hb)$ on $D_r\defin\{|\hb|<r\}$ such that, on $\AA_r$, $C_x(\hb)={}^t\Sigma'_x(\hb)\ovv{\Sigma''_x(\hb)}$. Therefore, the proof consists in finding $\Sigma',\Sigma''$ depending in a $C^\infty$ way on~$x$. The conclusion follows by checking that the formula given in \loccit for $\Sigma'$ depends in a $C^\infty$ way on $x$.
\end{proof}

\subsection{Variation of twistor structures}
The notion of variation of twistor structure was introduced by C\ptbl Simpson \cite{Simpson97}. He emphasized that it is convenient to express the relation between $D$ and $\VD$ in the definition of a harmonic bundle by adding a parameter $\hb$ and to express it with the notion of a $\hb$-connection, already suggested by P\ptbl Deligne.

By a real-analytic (\resp $C^\infty$) \emph{variation of twistor structure} on $M$ we mean the data of a triple $(\cH',\cH'',\cCS)$ defining a real-analytic (\resp $C^\infty$) family of holomorphic bundles on $\PP^1$ as above, such that each of the holomorphic bundles $\cH',\cH''$ is equipped with a \emph{relative} holomorphic connection
\[
\nabla:\cH'('')\to\frac1\hb\,\Omega^1_{M\times\Omega_0/\Omega_0}\otimes_{\cO_{M\times\Omega_0}}\cH'('')
\]
which has a pole along $\hb=0$ and is integrable. Moreover, the pairing $\cCS$ has to be compatible (in the usual sense) with the connections, \ie
\[
d'\cCS(m',\ovv{m''})=\cCS(\nabla m',\ovv{m''})\quad\text{and}\quad d''\cCS(m',\ovv{m''})=\cCS(m',\ovv{\nabla m''}).
\]
Let us note that we can define $\ovv\nabla$ as
\[
\ovv\nabla:\ovv{\cH''}\to \hb\Omega^1_{\ov M\times\Omega_\infty/\Omega_\infty}\otimes_{\cO_{\ov M\times\Omega_\infty}}\ovv{\cH''}.
\]
If we regard $\cCS$ as a $C^{*,\an}$-linear isomorphism
\begin{equation}\label{eq:cC}
\cCS:\cC^{*,\an}_{M_\RR\times\bS}\otimes_{\cO_{\ov M\times\bS}}\ovv{\cHS''}\isom\cC^{*,\an}_{M_\RR\times\bS}\otimes_{\cO_{M\times\bS}}\cHS^{\prime\vee},
\end{equation}
the compatibility with $\nabla$ means that $\cCS$ is compatible with the connection $d'+\ovv\nabla$ on the left-hand term and $\nabla^\vee+d''$ on the right-hand term, where $d',d''$ are the standard differentials with respect to $M$ only. The adjoint $(\cH',\cH'',\cCS)^\dag$ is defined as $(\cH'',\cH',\cCS^\dag)$, with
\[
\cCS^\dag(m'',\ovv{m'})\defin\ovv{\cCS(m',\ovv{m''})}.
\]
With respect to \eqref{eq:cC}, we can write $\cCS^\dag=\ovv\cCS{}^\vee$.

We say that the variation is
\begin{itemize}
\item
\emph{Hermitian} if $\cH''=\cH'$ and $\cCS$ is ``Hermitian'', \ie $\cCS^\dag=\cCS$;
\item
\emph{pure of weight~$0$} if the restriction to each $x\in M$ defines a \emph{trivial} holomorphic bundle on $\PP^1$.
\end{itemize}

\begin{lemme}[C\ptbl Simpson \cite{Simpson97}]\label{lem:simpson}
We have an equivalence between variations of Hermitian pure twistor structures of weight~$0$ and harmonic Higgs bundles, by taking $\PP^1$-global sections.
\end{lemme}

\begin{proof}[Sketch of proof ({\normalfont\itshape see \cite{Simpson97} for details, \cf also {\cite[Lemma 2.2.2]{Bibi01c}}})]\mbox{}
Let us start with a variation of Hermitian pure twistor structure of weight~$0$, that we denote by $(\cH',\nabla,\cCS)$. According to Proposition \ref{prop:familyP1}, the $C^{*,\an}$ bundle $\cH$ it defines is the pull-back by $p$ of some $C^*$ bundle, whose \emph{conjugate} bundle we now denote by $H$. By definition, $\ov H$ is naturally included in the $C^{*,\an}$ bundles associated to $\cH^{\prime\vee}$ and $\ovv{\cH'}$, and the inclusions are compatible with $\cCS$. Taking conjugates, $H$ is naturally included in the $C^{*,\an}$ bundles associated to $\ovv{\cH^{\prime\vee}}$ and $\cH'$. The natural duality pairing $\cH'\otimes\cH^{\prime\vee}\to\cO_{M\times\Omega_0}$ induces a pairing $h:H\otimes_{\cC^*_{M_\RR}}\ov H\to\cC^*_{M_\RR}$. That $\cCS$ is Hermitian implies that $h$ is Hermitian in the usual sense.

One defines $E$ as the restriction of $\cH'$ at $\hb=0$, $\Phi$ as the residue of $\nabla$ along $\hb=0$, $D+\Phi+\Phi^\dag$ as the restriction to $\hb=1$ of $\nabla+d''$.

Conversely, let us be given a harmonic Higgs bundle $(E,h,\Phi)$. Let $H$ be the $C^\infty$-bundle associated with $E$ and let $\cH_0=\cC^{*,\an}_{M_\RR\times\Omega_0}\otimes_{p^{-1}\cC^*_{M_\RR}}H$ be the pull-back of $H$ by $p:M\times\Omega_0\to M$ (where, as above, $*$ is for $\Ran$ or $\infty$, depending whether $h$ is real analytic or $C^\infty$). This bundle is equipped with the $d''$-operator $d''+\hb\Phi^\dag$, which defines a holomorphic bundle $\cH'=\ker (d''+\hb\Phi^\dag)$, equipped with the meromorphic connection $\nabla=(D'+\Phi/\hb)_{|\cH'}$. If we regard $h$ as a $\cC^*_{M_\RR}$-linear morphism $H\otimes_{\cC^*_{M_\RR}}\ov H\to\cC^*_{M_\RR}$, then we can extend $h$ as
\begin{equation}\label{eq:hS}
h_\bS:\cH_{0|M_\RR\times\bS}\otimes_{\cC^{*,\an}_{M_\RR\times\bS}}\ovv{\cH_{0|M_\RR\times\bS}}\to\cC^{*,\an}_{M_\RR\times\bS}
\end{equation}
by $\cC^{*,\an}_{M_\RR\times\bS}$-linearity and define $\cCS$ to be its restriction to $\cHS'\otimes_{\cO_\bS}\ovv{\cHS'}$.
\end{proof}

\begin{definition}
Let $(\cH',\nabla,\cCS)$ be a variation of Hermitian pure twistor structure of weight~$0$. We say that it is a \emph{variation of polarized pure twistor structure of weight~$0$} if $h$ is positive definite.
\end{definition}

\subsection{Supplementary structures on a variation of twistor structure}\label{subsec:integrability}\mbox{}
We introduce new structures in a way analogous to that of \S\ref{subsec:supplstruct}. We assume that $(\cH',\nabla,\cCS)$ is a variation of Hermitian pure twistor structure of weight~$0$.

\subsubsection*{Real structure}
Let $\cP$ be a nondegenerate $\cO_{M\times\Omega_0}$-linear morphism $\cH'\otimes_{\cO_{M\times\Omega_0}}\iota^*\cH'\to\cO_{M\times\Omega_0}$, or in other words an isomorphism $\cH^{\prime\vee}\isom\iota^*\cH'$. Let us notice that, on $\bS$, $\sigma$ and $\iota$ coincide. Therefore, restricting $\cP$ to $M\times\bS$ and composing with $\cCS$ (\cf \eqref{eq:cC}) we get an isomorphism
\begin{equation}\label{eq:cK}
\ov\cKS:\cC_{M_\RR\times\bS}^{*,\Ran}\otimes_{\cO_{\ov M\times\bS}}\ov{\cHS'}\isom\cC_{M_\RR\times\bS}^{*,\Ran}\otimes_{\cO_{M\times\bS}}\cHS'.
\end{equation}
(Let us note that we use the usual conjugation functor here.)

We say that $\cP$ is a \emph{real structure} if the following holds:
\begin{enumerate}
\item
$\cKS\ov\cKS=\id$,
\item
$\cP$ is $\iota$-Hermitian, \ie $\cP=\iota^*\cP^\vee$,
\item
$\cP$ is compatible with $\nabla$ and $\iota^*\nabla^\vee$.
\end{enumerate}

\begin{remarque}\label{rem:cP}
We could have defined the real structure directly from $\cKS$, but we emphasize here that the corresponding $\cPS$ should be the restriction of a $\iota$\nobreakdash-sesquilinear form $\cP$ on $\cH'$. In the integrable case below, $\cPS$ defines $\cP$ on $M\times\Omega_0^*$, so our assumption here means that we are given an extension of $\cP_{|\Omega_0^*}$ at $\hb=0$.
\end{remarque}

\begin{lemme}\label{lem:simpsonreel}
The equivalence of Lemma \ref{lem:simpson} specializes to objects with real structures.\qed
\end{lemme}

\subsubsection*{Integrability}
A variation of twistor structure $(\cH',\cH'',\cCS)$ is \emph{integrable} if the relative connection $\nabla$ on $\cH'$ comes from an absolute connection also denoted by $\nabla$, which has Poincar\'e rank one (\cf \cite[Chap\ptbl7]{Bibi01c}). In other words, $\hb\nabla$ should be an integrable meromorphic connection on $\cH',\cH''$ with a logarithmic pole along $\hb=0$. We also ask for a supplementary compatibility property of the absolute connection with the pairing in the following way:
\[
\hb\nabla_{\partial_\hb}\cCS(m',\ovv{m''})=\cCS(\hb\nabla_{\partial_\hb}m',\ovv{m''})-\cCS(m',\ovv{\hb\nabla_{\partial_\hb}m''}).
\]

\begin{exemple}\label{exem:PHS}
Assume $M$ is a point. A Hermitian pure twistor structure of weight~$0$ corresponds, according to Lemma \ref{lem:simpson}, to a vector space $H$ with a nondegenerate Hermitian form $h$. The integrability condition consists in giving a meromorphic connection on the trivial vector bundle $\cO_{\PP^1}\otimes H$ which has a pole of order at most~$2$ at~$0$ and~$\infty$ and which satisfies a compatibility condition with~$h$. The connection can be written as $d+(U/\hb-Q-\hb U')d\hb/\hb$, where $U,U'$ and $Q$ are endomorphisms of~$H$. Compatibility with the pairing $\cCS$ now translates as $U'=U^\dag$ and $Q=Q^\dag$, where~$^\dag$ means adjoint with respect to $h$. Assume that $h$ is positive definite. Then $Q$ is semi-simple and its eigenvalues are real, so that $H$ decomposes with respect to the eigenvalues of $Q$ as $\bigoplus_{\alpha\in[0,1[}\bigoplus_{p\in\ZZ}H_{\alpha+p}$. If we set $H^{p,-p}=\bigoplus_{\alpha\in[0,1[}H_{\alpha+p}$, we get a polarized complex Hodge structure of weight~$0$ on~$H$. The endomorphism $Q$ induces a semi-simple endomorphism on each $H^{p,-p}$ with eigenvalues in $[p,p+1[$. We are also left with an endomorphism $U$ of $H$ with no particular relation with the polarized Hodge structure.

In case we start with a polarized complex Hodge structure of weight~$0$ on $H$, we choose $U=0$ and $Q=\bigoplus_{p\in\ZZ}p\id_{H^{p,-p}}$ to recover an integrable polarized pure twistor structure of weight~$0$ (\cf Example \ref{exem:VHS}).
\end{exemple}

Lemma \ref{lem:simpson} is now made precise in the case of integrable twistor structures.

\begin{lemme}[C\ptbl Hertling \cite{Hertling01}, \cf also {\cite[Cor\ptbl7.2.6]{Bibi01c}}]\label{lem:hertling}
The equivalence of Lemma \ref{lem:simpson} specializes to an equivalence between integrable variations of Hermitian pure twistor structures of weight~$0$ and integrable harmonic Higgs bundles, \ie harmonic Higgs bundles equipped with endomorphisms $\cU,\cQ$ satisfying \eqref{eq:higgsint1}--\eqref{eq:higgsint5}.\qed
\end{lemme}

\begin{remarque}
The restriction to each point of $M$ of an integrable variation of polarized pure twistor structures of weight~$0$ is an integrable polarized pure twistor structure of weight~$0$. Therefore, it is in particular a polarized complex Hodge structure plus supplementary data (Example \ref{exem:PHS}). Nevertheless, it may not correspond, in general, to a variation of polarized Hodge structure of weight~$0$.

A sufficient condition to get such a variation of Hodge structure is that $\cU\equiv0$: in such a case, according to \eqref{eq:higgsint2} and \eqref{eq:higgsint5}, $\cQ$ is annihilated by $D$ and its eigenvalues are constant (use an orthonormal frame for $h$ in which the matrix of $\cQ$ is diagonal); the eigenspace decomposition of $H$ is left stable by $D$ and is shifted by one by $\Phi$, according to \eqref{eq:higgsint4}.

Such a condition is satisfied if we assume that $M$ is compact K\"ahler, so that the notion of \emph{integrable} variation of polarized twistor structure of weight~$0$ is a special case of that of variation of polarized complex Hodge structure (\cf \cite[Cor\ptbl7.2.8]{Bibi01c}).

The same holds if $M$ is a punctured compact Riemann surface and the behaviour at the punctures is \emph{tame} (\cf \cite{Simpson90} for this notion): this has been proved by C\ptbl Hertling and Ch\ptbl Sevenheck.\footnote{This was obtained in \cite[Cor\ptbl 7.2.10]{Bibi01c} as a consequence of the conjecture \cite[Conj\ptbl 7.2.9]{Bibi01c}. Ch\ptbl Sevenheck sent a proof of this conjecture due to C\ptbl Hertling and himself to the author.}

One can expect that such a result holds in the higher dimensional case, if we keep the tameness assumption.

Therefore, when $M$ is quasi-projective, the notion of \emph{integrable} variation of polarized twistor structure of weight~$0$ should be a new object only when a nontame (\ie wild) behaviour occurs at infinity on $M$. In other words, this notion can be regarded as the right replacement of that of a polarized complex variation of Hodge structure of weight~$0$ when one wants to take into account a wild behaviour at infinity.
\end{remarque}

Let us now come back to the general notion of an integrable Hermitian variation of twistor structure. One can remark that, as $\nabla$ is integrable, the restriction of $(\cH',\nabla)$ to $M\times\Omega_0^*$ (where $\Omega_0^*\defin\Omega_0\moins\{0\}$) is determined up to isomorphism by giving the corresponding locally constant sheaf $\ker\nabla$. Similarly, on $M_\RR\times\bS$, the pairing $\cCS$ is determined by its restriction to the corresponding local systems. Therefore, giving an integrable variation of twistor structure with a nondegenerate Hermitian pairing is equivalent to giving
\begin{enumerate}
\item
a holomorphic vector bundle $\cH'$ on $M\times\Omega_0$ with a meromorphic connection~$\nabla$ having Poincar\'e rank one along $M\times\{0\}$ (\ie $\hb\nabla$ has at most a logarithmic pole along $M\times\{0\}$),
\item
a nondegenerate pairing $\cCS:\ker\nabla_{|M\times\bS}\otimes_\CC\sigma^{-1} \ov{\ker\nabla_{|M\times\bS}}\to\CC_{M\times\bS}$.
\end{enumerate}
The second part of the data is of a purely topological nature, as $\ker\nabla_{|M\times\Omega_0^*}$ is a locally constant sheaf. Moreover, recall that $\iota$ coincides with $\sigma$ on $\bS$. As a direct application of Proposition \ref{prop:familyP1} we get:

\begin{corollaire}[C\ptbl Hertling \cite{Hertling01}]\label{cor:deformation}
Let us assume that $M$ is $1$-connected and let \hbox{$x^o\in M$}. Let $(\cH',\nabla)$ be a holomorphic bundle on $M\times\Omega_0$ with a meromorphic connection~$\nabla$ having Poincar\'e rank one along $M\times\{0\}$, set $\cLS=\ker\nabla_{|M\times\bS}$, $\cLS^o=\cL_{|\{x^o\}\times \bS}$ and let \hbox{$\cCS^o:\cLS^o\otimes_\CC\iota^{-1}\ov{\cLS^o}\to\CC_\bS$} be a nondegenerate Hermitian pairing. Let \hbox{$\cCS:\cLS\otimes_\CC\iota^{-1}\ov{\cLS}\to\CC_{M\times\bS}$} be the unique nondegenerate Hermitian pairing extending~$\cCS^o$. It defines a real-analytic family of holomorphic bundles parametrized by~$M$. Let us moreover assume that the twistor structure at~$x^o$ corresponding to these data is Hermitian and pure of weight~$0$.

Then there exists a (possibly empty) real analytic subvariety $\Theta\not\ni x^o$ of $M$ such that, on the connected component of $M\moins\Theta$ containing~$x^o$, the Hermitian variation of twistor structure $(\cH',\nabla,\cCS)$ is pure of weight~$0$. If moreover it is polarized at~$x^o$, it is polarized on this connected component.\qed
\end{corollaire}

\begin{remarques}\label{rem:deformation}\mbox{}
\begin{enumerate}
\item
This corollary shows that, given $(\cH',\nabla)$ having Poincar\'e rank one along \hbox{$\hb=0$}, it underlies an integrable variation of polarized pure twistor structure in some neighbourhood of $x^o$ if and only if the fibre at $x^o$ can be endowed with the structure of an integrable polarized twistor structure of weight~$0$.
\item
Let us start with an integrable pure twistor structure $(\cH^{\prime o},\nabla^o,\cCS^o)$ of weight $0$ and let us assume that the residue of $\hb\nabla^o$ at $\hb=0$ is \emph{regular semi-simple} (\ie has pairwise distinct eigenvalues). Then (\cf \cite{Malgrange83db}, see also \cite[Th\ptbl III.2.10]{Bibi00}) there exists a universal deformation $(\cH',\nabla)$ of $(\cH^{\prime o},\nabla^o)$ parametrized by the universal covering $M$ of $\CC^d\moins\text{diagonals}$, if $d$ is the rank of $\cH^{\prime o}$, and $\cCS^o$ extends in a unique way to a Hermitian pairing $\cCS$. We can apply the conclusion of Corollary \ref{cor:deformation} to the object $(\cH',\nabla,\cCS)$.
\item
Let us denote by $L^o$ the space of global multivalued sections of $\cL^o$ (\ie of global multivalued $\nabla^o$-horizontal sections of $\cH^{\prime o}_{|\{x^o\}\times\CC^*}$). Giving $\cCS^o$ is equivalent (after choosing the lifting $\zeta\mto\zeta+\pi i$ of $\iota:\hb=e^\zeta\mto-\hb$) to giving a nondegenerate Hermitian pairing $C^o:L^o\otimes_\CC\ov{L^o}\to\CC$. However, the purity condition and the positivity condition are not conditions on $C^o$ (\ie on $\nabla^o$-horizontal sections). They much depend on the holomorphic bundle $\cH^{\prime o}$. More precisely, the purity and positivity conditions are equivalent to the existence of a global holomorphic frame $\varepsilong^o$ of $\cH^{\prime o}$ which is orthonormal with respect to $\cCS^o:\cH^{\prime o}_{|\bS}\otimes_{\cO_\bS}\ovv{\cH^{\prime o}_{|\bS}}\to\cO_\bS$. In other words, for any $\hbo\in\bS$, the matrix $\big(\cCS^o(\epsilon_i^{o,(\hbo)},\epsilon_j^{o,(-\hbo)})\big)$ should be equal to the identity matrix (\cf \cite[Rem\ptbl 2.1.4 and 2.2.3]{Bibi01c}).
\end{enumerate}
\end{remarques}

\subsubsection*{Integrability with real structure}
Let $(\cH',\nabla,\cCS,\cP)$ be a variation of Hermitian pure twistor structure of weight~$0$ with real structure $\cP$. We say that it is \emph{integrable} if $(\cH',\nabla,\cCS)$ is integrable, as defined above, and $\cP$ is compatible with the absolute connections (not only the relative ones).

The conjunction of Lemmas \ref{lem:simpsonreel} and \ref{lem:hertling} holds. Similarly, starting with an \emph{integrable} $(\cH',\nabla,\cP)$, we get as in corollary \ref{cor:deformation} that, if $M$ is $1$-connected, an integrable variation of real Hermitian twistor structure on $M$ with underlying $(\cH',\nabla,\cP)$ is determined by a real structure $\cKS^o:\cLS^o\isom\ov{\cLS^o}$ at $x^o$, and if it is pure of weight $0$ (\resp pure of weight $0$ and polarized) at $x^o$, it is so on $M\moins\Theta$.

\subsubsection*{Potentiality}
We now consider the analogue of the potentiality property. Let $(\cH',\nabla,\cCS)$ be a variation of Hermitian pure twistor structure of weight~$0$. By a \emph{potential} we will mean an isomorphism of holomorphic bundles $p^*E\isom\cH'$ inducing the identity when restricted to $\hb=0$. Recall that $\cH'$ is the holomorphic bundle $(p^*H,d''+\hb\Phi^\dag)$. On the other hand, $p^*E=(p^*H,d'')$. Therefore, such an isomorphism can be expanded as $\id_H+\hb B_1+\cdots+\hb^kB_k+\cdots$, where the $B_k$ are $C^\infty$ endomorphisms of $H$ and the condition is
\[
(\id_H+\hb B_1+\cdots)^{-1}\circ d''\circ (\id_H+\hb B_1+\cdots)=d''+\hb\Phi^\dag.
\]
In particular, the order-one condition is $d''B_1=\Phi^\dag$, that is, the potentiality condition for the harmonic Higgs bundle corresponding to $(\cH',\nabla,\cCS)$. But the potentiality condition for the variation of twistor structure is stronger than that for harmonic Higgs bundles.

\subsection{A criterion for a Saito structure to be harmonic}
Let us consider a Saito structure $(E,\nablaf,g,\Phi,\cU,\cV)$ of weight $0$. We consider on $M\times\Omega_0$ the pull-back bundle $p^*E$ via $p:M\times\Omega_0\to M$, equipped with the meromorphic connection
\begin{equation}\label{eq:nablafnabla}
\nabla=p^*\nablaf+\frac{\Phi}{\hb}+\Big(\frac{\cU}{\hb}-\cV\Big)\frac{d\hb}{\hb}.
\end{equation}
The bilinear form $g$ defines then a nondegenerate bilinear pairing
\begin{equation}\label{eq:cP}
\cP:(p^*E,\nabla)\otimes_{\cO_{M\times\Omega_0}}\iota^*(p^*E,\nabla)\to(\cO_{M\times\Omega_0},d),
\end{equation}
which is Hermitian with respect to the involution $\iota$. If we denote by $\cL$ the local system $\ker\nabla_{|M\times\CC^*}$ and by $\cLS$ its restriction to $M\times\bS$, where $\bS=\{\hb\mid |\hb|=1\}$, we get a nondegenerate pairing $\cPS:\cLS\otimes_{\CC}\iota^{-1}\cLS\to\CC_{M\times\bS}$.

\begin{corollaire}[{\cf \cite[Th\ptbl5.15]{Hertling01}}]\label{cor:critharmonic}
In such a situation, assume moreover that there exists a real structure $\cKS$ on $\cLS$, that is, an isomorphism $\cKS:\cLS\to\ov{\cLS}$ such that $\cKS\ov\cKS=\id$, such that, setting $\cCS(\cbbullet,\cKS\cbbullet)=\cP(\cbbullet,\cbbullet)$, $(p^*E,\nabla,\cCS)$ defines an integrable variation of Hermitian pure twistor structure of weight~$0$. Then the Saito structure $(E,\nablaf,g,\Phi,\cU,\cV)$ is harmonic.
\end{corollaire}

\begin{proof}
As the variation of twistor structure is defined on $p^*E$, we have a given identification $\cH'=p^*E$ and the potentiality property is fulfilled. The corollary is then a consequence of the extensions of Lemma \ref{lem:simpson} given above.
\end{proof}

\begin{remarques}\mbox{}
\begin{enumerate}
\item
By the corollaries \ref{cor:deformation} and \ref{cor:critharmonic}, if there exists $x^o\in M$ such that the data $(\cO_{\Omega_0}\otimes_\CC E^o,d+(\cU^o/\hb-\cV)d\hb/\hb,\cP^o,\cKS^o)$ defines a pure twistor structure of weight $0$, then $(E,\nablaf,g,\Phi,\cU,\cV)$ is harmonic on some neighbourhood of $x^o$. The same conclusion holds with the polarization property.
\item
This proposition can be applied to the Saito structure attached to a Frobenius manifold, and gives a criterion for a Frobenius manifold to be harmonic.
\end{enumerate}
\end{remarques}

\section{Twistor structures through Fourier-Laplace transform}\label{sec:Fourier-Laplace}

According to Corollary \ref{cor:deformation}, an important step in constructing an integrable variation of pure twistor structure of weight~$0$ is to construct it at one point, \ie to construct an integrable pure twistor structure of weight~$0$. Checking that $(\cH^{\prime o},\nabla^o,\cCS^o)$ is such an object is not straightforward. We indicate in this section how to produce such an integrable pure twistor structure of weight~$0$ starting from a variation of polarized Hodge structure of weight~$0$ on a punctured Riemann sphere, by using Fourier-Laplace transform. This will produce an integrable polarized pure twistor structure of weight~$0$: the polarization property is essential here, to get the purity property.

\subsection{Laplace transform}\label{subsec:Laplace}
Let $V^\an$ be a holomorphic vector bundle on the punctured affine line $\Afu\moins P$, where $P=\{p_1,\dots,p_n\}$ is a finite set of points. We assume that $V^\an$ is equipped with a holomorphic connection $\nabla^\an$. By a classical result of Deligne, there exists a unique meromorphic extension $(V,\nabla)$ of $(V^\an,\nabla^\an)$ on $\PP^1$ such that $\nabla$ has only regular singularities at $\{p_1,\dots,p_n,\infty\}$. Taking global sections of $V$ on $\PP^1$ gives a $\CC[t]$-module $M$, where $t$ denotes a chosen coordinate on $\Afu$. Moreover, the connection~$\nabla$ enables one to extend the $\CC[t]$ action to an action of the Weyl algebra $\CC[t]\langle\partial_t\rangle$ by setting $\partial_tm\defin\nabla_{\partial_t}m$.

The Laplace transform $\Fou M$ is the $\CC$-vector space $M$ equipped with the action of the Weyl algebra $\CC[\tau]\langle\partial_\tau\rangle$ defined by the formulas
\[
\tau m\defin\partial_tm,\quad\partial_\tau m=-tm.
\]
We will be mainly interested in the behaviour of $\Fou M$ near $\tau=\infty$ and we will set $\hb=1/\tau$, so that $\hb^2\partial_\hb=t$. It is known that $G\defin\CC[\tau,\tau^{-1}]\otimes_{\CC[\tau]}\Fou M$ is a free $\CC[\tau,\tau^{-1}]$-module, or equivalently a free $\CC[\hb,\hb^{-1}]$-module, \ie an algebraic vector bundle on the torus $\CC^*$. Moreover, the action of $t=-\partial_\tau=\hb^2\partial_\hb$ defines a connection $\Fou\nabla$ on it, which is known to have a regular singularity at $\tau=0$ and, in general, an irregular one at $\hb=0$ (see \eg \cite{Malgrange91} or \cite[Chap\ptbl V]{Bibi00}).

Assume now that we are given a $\CC[t]$-submodule of finite type $L\subset M$ such that $M=L+\partial_tL+\cdots+\partial_t^kL+\cdots$. Then we can regard $L$ as a subspace of $\Fou M$ and, taking localization, we can consider its image $L'$ in $G$. Let us then set $G_0^{(L)}=L'+\hb L'+\cdots+\hb^kL'+\cdots\subset G$. By construction, $G_0^{(L)}$ is a $\CC[\hb]$-submodule of~$G$ and, $L$ being a $\CC[t]$-module, it is naturally equipped with a compatible action of $\hb^2\partial_\hb$. One can show (see \eg \cite[\S V.2.c]{Bibi00}) that, because $\nabla$ has a regular singularity at infinity, $G_0^{(L)}$ is a $\CC[\hb]$-module of finite rank. As it is equipped with an action of $\hb^2\partial_\hb$ (because $L'$ is so), and as $G=\CC[\hb,\hb^{-1}]\otimes_{\CC[\hb]}G_0^{(L)}$, we find that the connection~$\Fou\nabla$ on $G$ has at most Poincar\'e rank one at $\hb=0$ (in any $\CC[\hb]$-basis of $G_0^{(L)}$, the matrix of $\Fou\nabla$ has a pole of order $\leq1$ because $G_0^{(L)}$ is stable by $\hb^2\partial_\hb$).

We will define $(\cH^{\prime o},\nabla^o)$ to be the analytization of $(G_0^{(L)},\Fou\nabla)$.

\subsection{Fourier transform of a sesquilinear pairing}\label{subsec:Fsesqui}
Let us assume that the holomorphic bundle $V^\an$ is equipped with a $\nabla$-\emph{horizontal} nondegenerate Hermitian pairing $\wt k$ (but not necessarily positive definite), that is, a $\cO_{\Afu\moins P}\otimes_\CC\cO_{\ov\Afu\moins P}$-linear morphism
\[
\wt k:V^\an\otimes_\CC\ov{V^\an}\to\cC^\infty_{\Afu\moins P}
\]
such that $d'\wt k(v_1,\ov{v_2})=\wt k(\nabla v_1,\ov{v_2})$ and $d''\wt k(v_1,\ov{v_2})=\wt k(v_1,\ov{\nabla v_2})$. Let us moreover assume that $\wt k$ extends to $M\otimes_\CC\ov M$ as a Hermitian sesquilinear pairing $k$ with values in the Schwartz space $\cS'(\Afu)$ of temperate distributions on $\Afu$ (\ie distributions on $\Afu$ having moderate growth at infinity), in a way compatible to the natural action of $\CC[t]\langle\partial_t\rangle\otimes_\CC\ov{\CC[t]\langle\partial_t\rangle}$ on $\cS'(\Afu)$.

\begin{remarque}[A criterion for $k$ to be nondegenerate]\label{rem:knondeg}
While the nondegeneracy property of $\wt k$ is easy to define (by restricting to points of $\Afu\moins P$), that of $k$ needs some care. The first step consists in considering the $\CC[t]\langle\partial_t\rangle$-module $M^\dag\defin\Hom_{\CC[\ov t]\langle\ov \partial_t\rangle}\big(\ov M,\cS'(\Afu)\big)$, where the $\CC[t]\langle\partial_t\rangle$-structure comes from that on $\cS'(\Afu)$. It is called the \emph{Hermitian dual} of $M$ and was introduced by M\ptbl Kashiwara in \cite{Kashiwara86}. Then giving $k$ is equivalent to giving a $\CC[t]\langle\partial_t\rangle$-linear morphism $M\to M^\dag$. We say that $k$ is \emph{nondegenerate} if this morphism is an isomorphism.

In general, if $M$ is holonomic, then $M^\dag$ is so: this is proved in \cite{Kashiwara86} if $M$ has only regular singularities and in \cite{Bibi99} in general (but on a manifold of dimension one, $\Afu$ here).

Moreover, the previous construction can be localized in the analytic topology of~$\PP^1$; in particular, we replace $\cS'(\Afu)$ with the sheaf of distributions on $\Afu$ having moderate growth at infinity. Then, in restriction to $(\Afu\moins P)^\an$, the nondegeneracy of~$k$ in the sense of $\cD$-modules is equivalent to that of $\wt k$ in the usual sense. In other words, saying that $\wt k$ is nondegenerate is equivalent to saying that the kernel and the cokernel of $k:M\to M^\dag$ are supported on $P$.

At this point, let us recall a definition: we say that $M$ is a \emph{minimal extension} (one also says an intermediate extension, or a middle extension, with reference to perverse sheaves) if $M$ has neither a nonzero submodule nor a nonzero quotient module supported on $P$.

As the construction $^\dag$ is functorial (in a contravariant way) for $\CC[t]\langle\partial_t\rangle$-modules (see \loccit), $M$ is a minimal extension iff $M^\dag$ is so. In such a case, the kernel and cokernel of $k$, if supported at $P$, must be zero. In other words, \emph{if $\wt k$ is nondegenerate and is extended to a sesquilinear pairing $k$ on the minimal extension $M$ of $(V,\nabla)$, then $k$ is also nondegenerate}.

This motivates the change of the definition of $M$ in \S\ref{subsec:Laplace}, and the replacement of the $\CC[t]\langle\partial_t\rangle$-module naturally associated to $M$ with its minimal extension.
\end{remarque}

\subsubsection*{Fourier transform of $k$}
We then denote by $\Fou k$ the composition of $k$ with the Fourier transform of temperate distributions with kernel $\exp(\ov{t\tau}-t\tau)\itwopi dt\wedge d\ov t$. By standard properties of the Fourier transform, $\Fou k$ is a $\CC[\tau]\langle\partial_\tau\rangle\otimes_\CC\ov{\CC[\tau]\langle\partial_\tau\rangle}$-linear morphism on $\Fou M\otimes_\CC\iota^*\ov{\Fou M}$ (note the $\iota^*$ here, due to $\ov{\exp(\ov{t\tau}-t\tau)}=\exp-(\ov{t\tau}-t\tau)$).

Considering the Hermitian duality functor, we see that, composing with the Fourier transform $\cS'(\Afu_t)\to\cS'(\Afu_\tau)$ gives an isomorphism of $\CC$-vector spaces
\[
\Hom_{\CC[\ov t]\langle\ov \partial_t\rangle}\big(\ov M,\cS'(\Afu_t)\big)\to\Hom_{\CC[\ov\tau]\langle\ov \partial_\tau\rangle}\big(\Fou{\ov M},\cS'(\Afu_\tau)\big)
\]
which is in fact a canonical identification $\iota^*(\Fou M)^\dag=\Fou(M^\dag)$.

If we regard $k$ as a linear morphism $M\to M^\dag$, taking its Laplace transform gives a linear morphism $\Fou M\to \Fou(M^\dag)$, and using Fourier identification above gives $\Fou k:\Fou M\to\iota^*(\Fou M)^\dag$. In particular, if $k$ is nondegenerate, then so is $\Fou k$.

\subsubsection*{Twistorization}
In particular, restricting to $\tau\neq0$ and sheafifying in the analytic topology gives a $\nabla^o$-flat nondegenerate Hermitian pairing (with respect to the composition of the involution $\iota$ and the conjugation) $\cH^{\prime o}_{|\CC^*}\otimes_\CC\iota^*\ov{\cH^{\prime o}_{|\CC^*}}\to\cC^\infty_{\CC^*}$. By flatness, giving such a pairing is equivalent to giving the pairing induced on $\nabla^o$-horizontal sections
\[
\Fou k:\cL^o\otimes_\CC\iota^{-1}\ov{\cL^o}\to\CC_{\CC^*}.
\]
Last, as $\bS\hto\CC^*$ is a homotopy equivalence, giving this pairing is equivalent to giving the restricted pairing
\[
\cCS^o\defin\Fou k_{|\bS}:\cLS^o\otimes_\CC\iota^{-1}\ov{\cLS^o}\to\CC_\bS.
\]

In conclusion (and extending slightly the previous construction), \emph{given $(M,L,k)$, where $M$ is a holonomic $\CC[t]\langle\partial_t\rangle$-module having a regular singularity at infinity and which is a minimal extension at~$P$, $L$ is a $\CC[t]$-submodule of finite type of $M$ which generates it over $\CC[t]\langle\partial_t\rangle$, and $k$ is a nondegenerate sesquilinear Hermitian pairing on $M$ with values in $\cS'(\Afu)$, we obtain, by Fourier-Laplace transform, a twistor structure $(\cH^{\prime o},\nabla^o,\cCS^o)$}.

\begin{remarque}[Conjugate to Remark \ref{rem:cP}]
In the previous construction, $\cCS^o$ is the restriction to $\bS$ of a sesquilinear pairing $\Fou k$ which is defined at $\hb=0$ and taking value there in the space of distributions having a moderate growth $\hb=0$. More generally, given a free $\CC\{\hb\}$-module $\cG_0$ equipped with a connection having Poincar\'e rank one, and with a Hermitian nondegenerate sesquilinear pairing $c$ on $\cG_0\otimes_\CC\iota^*\ov{\cG_0}$ taking values in the space of germs of distributions with moderate growth at $\hb=0$, $c$ being compatible with the $\partial_\hb$- and $\ov\partial_\hb$-actions, one can ask for a criterion ensuring that, if we extend flatly $\cG_0$ and $c$ to $\Omega_0$, the triple $(\cG_0,\partial_\hb,c_\bS)$ is pure of weight $0$ (and polarized).

Let us also recall that $\cG_0[1/\hb]$ is completely determined in terms of the formal meromorphic connection $\wh\cG_0[1/\hb]\defin\CC\lcr\hb\rcr\otimes_{\CC\{\hb\}}\cG_0[1/\hb]$ and the Stokes structure, and that $\wh\cG_0[1/\hb]$ has a decomposition into simple objects (Turrittin-Levelt decomposition). Moreover, according to a result of Malgrange, $\cG_0$ is completely determined by $\wh\cG_0\defin\CC\lcr\hb\rcr\otimes_{\CC\{\hb\}}\cG_0$.

One can ask whether $c$ can be determined by the data of sesquilinear pairings on the objects of the Turrittin-Levelt decomposition of $\wh\cG_0[1/\hb]$ (that we could denote by $\wh c$) and a correspondence between Stokes matrices and their transpose-conjugate (to go from $c$ to a sesquilinear pairing on the objects of the Turrittin-Levelt decomposition, one can use \cite{Bibi99} after twisting by suitable exponential factors; on the other hand, describing the correspondence between the Stokes cocycle for $\cG_0[1/\hb]$ and $\cG_0[1/\hb]^\dag$ is the contents of the proof of \cite[Th\ptbl II.3.1.2]{Bibi99}). Would this be the case, it would be interesting to give a criterion for the purity or polarization of $(\cG_0,c_\bS)$ in terms of $(\wh\cG_0,\wh c)$ and the Stokes matrices (\cf \cite[Conjecture 10.2]{H-S06} for a precise conjecture in this direction and \S8 of \loccit for a more precise analysis).
\end{remarque}

\subsection{Application to variations of Hodge structures}

Let us assume that $(V^\an,\nabla)$ underlies a variation of polarized complex Hodge structure of weight~$0$ with flat Hermitian pairing $\wt k$ and associated metric~$h$ (\cf Example \ref{exem:VHS}). Let $\wt V$ be the subsheaf of $j_*V^\an$ ($j:\Afu\moins P\hto \PP^1$) consisting of sections whose $h$-norm has moderate growth near $P\cup\{\infty\}$, and let $\wt L\subset \wt V$ be the subsheaf consisting of sections of $\wt V$ whose $h$-norm is $L^1_\loc$ near any $p_i\in P$ (with respect to the local Euclidean metric at $p_i$). Classical results of Schmid's theory \cite{Schmid73} on variations of Hodge structures, suitably extended to variations of complex Hodge structures, imply that $\wt V$ is an algebraic vector bundle on $\Afu\moins P$ with a connection having only regular singularities, and $\wt L$ is an algebraic vector bundle on $\Afu$.

Last, let us denote by $M$ the $\CC[t]\langle\partial_t\rangle$-submodule of $V\defin\Gamma(\PP^1,\wt V)$ generated by $L=\Gamma(\PP^1,\wt L)$. Then it can be shown that $M$ is a minimal extension at each point of $P$.

Moreover, the sesquilinear pairing $\wt k$ extends to $M$ with values in $\cS'(\Afu)$ (see the details in \cite[\S3g]{Bibi05}), in a nondegenerate way as we have seen in Remark \ref{rem:knondeg}.

\begin{theoreme}[\cite{Bibi05}]\label{th:FourierHodge}
Under these conditions, the associated twistor structure $(\cH^{\prime o},\nabla^o,\cCS^o)$ obtained by Fourier-Laplace transform and twistorization is pure of weight~$0$ and polarized.
\end{theoreme}

\begin{proof}[Sketch of proof]
Let us first explain in other terms the construction above, in particular the twistorization step, which looks mysterious. To the variation of Hodge structure one associates an integrable variation of twistor structure, by adding the external parameter $\hb$ (\cf Lemma \ref{lem:hertling}).

There is a well-defined notion of Fourier-Laplace transform of a variation of twistor structure (integrable or not) and we are mainly interested in the fibre at $\tau=1$ of this Fourier-Laplace transform (where $\tau$ is the parameter introduced in \S\ref{subsec:Laplace}). The fundamental result, independent of the integrability condition, is

\begin{theoreme}\label{th:FLtw}
The Fourier-Laplace transform of a tame variation of polarized pure twistor structure of weight $0$ is a variation of polarized pure twistor structure of weight $0$ on $\CC^*=\{\tau\neq0\}$ with a tame behaviour at $\tau=0$. In particular, the fibre at $\tau=1$ of this variation is a polarized pure twistor structure of weight $0$.
\end{theoreme}

On the other hand, the behaviour at $\tau=\infty$ is usually not tame, and is understood in particular cases by S\ptbl Szabo \cite{Szabo04}.

Let us now use the integrability condition. Starting from an integrable variation of twistor structure, its Fourier-Laplace transform remains integrable. Therefore, the fibre at $\tau=1$ is an integrable twistor structure. When we more specifically start from a variation of Hodge structure, where we have a Hodge filtration, this integrable twistor structure is canonically identified with the twistor structure obtained from $(G_0^{(L)},\Fou\nabla,\Fou k)$. In other words, although the variation of twistor structure uses two a~priori distinct variables $\tau^{-1}$ and $\hb$, this variation is produced from the fibre at $\tau=1$ by rescaling the variable $\hb$, the rescaling parameter being $\tau$. This is analogous (without considering any real structure however) to the orbits of TERP structures considered in \cite{H-S06} and explains the rather mysterious twistorization construction above (\cf the proof of the main theorem 1.32 in \cite{Bibi05}).

Let us now come back to the proof of Theorem \ref{th:FLtw}, which is the main point for getting purity. Let us emphasize that we get purity together with polarization, and that, without the polarization assumption on the original variation of twistor structure, the $L^2$ techniques are not enough.

According to Lemma \ref{lem:simpson}, we consider a holomorphic bundle $E$ on $\Afu\moins P$, equipped with a Hermitian metric~$h$ and a holomorphic Higgs field $\Phi$. We assume that $(E,h,\Phi)$ is a \emph{tame harmonic bundle} in the sense of Simpson \cite{Simpson90} (\ie $(E,h,\Phi)$ is a harmonic Higgs bundle on $\Afu\moins P$ and we impose a growth condition for the eigenvalues of the Higgs field at the points of $P\cup\{\infty\}$).

Let us consider the exponentially twisted harmonic bundle $(E,h,\Phi-dt)$. This remains a harmonic bundle which is tame at $P$, but \emph{wild} at $\infty$. The main point in the proof of Theorem \ref{th:FLtw} is to find (using the notation of Theorem \ref{th:FourierHodge}) a finite dimensional $\CC$-vector space $H^o$ in $\cH^{\prime o}$ which contains an orthonormal basis for $\cCS^o$. This is done using $L^2$-theory:

\begin{theoreme}[\cite{Bibi04}]\label{th:L2FL}
The space of $L^2$ harmonic sections of $E\otimes\cA^1_{\Afu\moins P}$ ($1$-forms with values in $E$), with respect to the metric $h$ and a metric on $\Afu\moins P$ equivalent to the Poincar\'e metric near $P\cup\{\infty\}$, and with respect to the Laplace operator of $d''+\Phi-dt+\hb(D'+\Phi^\dag-d\ov t)$ is finite dimensional and independent of $\hb$.
\end{theoreme}

\begin{remarque}
The proof of this theorem that I first gave by using a degeneration argument $\tau\to0$ was only valid for $(E,h,\Phi-\tau dt)$ for $|\tau|$ small (depending on $(E,h,\Phi)$). The proof I~gave can be made correct (\cf the erratum in \cite{Bibi04}) by using in particular an argument whose proof is due to the the referee of \cite{Bibi05} (note also that S.~Szabo proves a similar result in \cite[Lemma 2.32]{Szabo04}, with different methods however):
\par\nopagebreak
\emph{The harmonic sections have an exponential decay when $t\to\infty$}.
\end{remarque}

Regarding the parameter $\hb$ in Theorem \ref{th:L2FL} as the external parameter of a twistor structure (the metric being given by the $L^2$ metric), the result of this theorem is the $L^2$-construction of a polarized pure twistor structure of weight $0$.

In order to realize our $L^2$-twistor structure as the fibre at $\tau=1$ of the Fourier-Laplace transform of the original variation of twistor structure, one uses a technique which goes back to S\ptbl Zucker \cite{Zucker79} which provides an identification with $L^2$ meromorphic sections (see the details in \cite{Bibi04,Bibi05}). In fact, this identification is already needed in the proof of Theorem \ref{th:L2FL} to ensure the finite dimension of the spaces of $L^2$ sections which are involved.
\end{proof}

\section{The harmonic Frobenius manifold attached to a Laurent polynomial}
\label{sec:laurent}
Let $U$ be the torus $\spec\CC[u_1,\dots,u_n,u_1^{-1},\dots,u_n^{-1}]\simeq(\CC^*)^n$ and let $f\in\CC[u_1,\dots,u_n,u_1^{-1},\dots,u_n^{-1}]$ be a Laurent polynomial. All along this section, we will assume that $f$ is convenient and nondegenerate in the sense of Kouchnirenko \cite{Kouchnirenko76}. Let us recall that the convenience assumption means that the origin belongs to the interior of the Newton polyhedron of $f$. A nice example is $f_w(u_1,\dots,u_n)=u_1+\cdots+u_n+1/u_1^{w_1}\cdots u_n^{w_n}$ with $w_i\in\NN^*$ (and, for simplicity, $\gcd(w_1,\dots,w_n)=1$).

We will also use a stronger convenience assumption.

\begin{definition}[A\ptbl Douai]
We say that $f$ is strongly convenient if, denoting by $0=\alpha_0,\alpha_1,\dots,\alpha_r$ the integral points contained in the interior of the Newton polyhedron of $f$, the polynomial ring on $u^{\alpha_1},\dots,u^{\alpha_r}$ surjects onto $\CC[u,u^{-1}]/(\partial f)$.
\end{definition}

The assumption holds for instance if there exists a $\ZZ$-basis $e_1,\dots,e_n$ of $\ZZ^n$ such that $\{\pm e_1,\dots,\pm e_n\}$ belong to the interior of the Newton polyhedron of $f$.

To any convenient and nondegenerate Laurent polynomial is attached in a canonical way a Frobenius manifold, under one of the following assumptions:

\refstepcounter{equation}\label{cond:semissimple}
\par\smallskip\noindent\eqref{cond:semissimple}\enspace
$f$ has only nondegenerate critical points and the critical values are distinct,
\refstepcounter{equation}\label{cond:strongconv}
\par\noindent\eqref{cond:strongconv}\enspace
$f$ is strongly convenient.

\smallskip
Let us note that the Laurent polynomial $f_w$ satisfies the first assumption but not the second one.

We recall in \S\ref{subsec:Douai} the construction of this canonical Frobenius manifold, which is due to A\ptbl Douai (\cite{Douai05,Douai07}) in the strongly convenient case \eqref{cond:strongconv}. In \S\ref{subsec:tame}, we enrich this structure to a structure of harmonic Frobenius manifold. The main argument consists in endowing the Brieskorn lattice of $f$ with a polarized pure twistor structure. This point is not specific to convenient and nondegenerate Laurent polynomials, but is a general result for cohomologically tame functions on affine manifolds. We thus take the general point of view in \S\ref{subsec:tame}.

\subsection{The Frobenius manifold attached to a Laurent polynomial}\label{subsec:Douai}

Let~$f$ be a convenient and nondegenerate Laurent polynomial. The $\CC$-vector space $\CC[u,u^{-1}]/(\partial f)$ is finite dimensional and its dimension $\mu(f)$ is the sum of the Milnor numbers of $f$ at each of its critical points. We will construct a canonical Frobenius manifold structure on the analytic germ $M$ along a submanifold $N$ of the $\CC$-vector space $\CC[u,u^{-1}]/(\partial f)$.

The method developed by A\ptbl Douai to construct such a structure consists in
\begin{enumerate}
\item
choosing a smooth affine subvariety $N\subset \CC[u,u^{-1}]/(\partial f)$ canonically attached to~$f$,
\item
constructing a Saito structure $(E,\nablaf,g,\Phi,\cU,\cV)$ on a holomorphic vector bundle $E$ of rank $\mu(f)$ on $N$ and canonically attached to $f$,
\item
giving a $\nablaf$-horizontal section $\omega$ of $E$, which is homogeneous with respect to $\cV$, such that the infinitesimal period map $\varphi_\omega:TN\to E$ is injective,
\item
apply the result of \cite{Malgrange83db} in Case \eqref{cond:semissimple} or \cite{H-M04} in Case \eqref{cond:strongconv} concerning the existence of a universal unfolding of a differential system to show that $\varphi_\omega$ extends in a unique way as a primitive homogenous section on some analytic neighbourhood $M$ of $N$ in $\CC[u,u^{-1}]/(\partial f)$,
\item
obtaining the Frobenius manifold structure on $M$ by carrying through $\varphi_\omega$ the Saito structure on $E$ to a Saito structure on $TM$.
\end{enumerate}

The important point is that, on $N$, the construction is completely algebraic, while the transcendental part of the construction is hidden in the application of the results of \cite{Malgrange83db} or \cite{H-M04}. This improves much the construction given in \cite{D-S02a}, where a universal unfolding of the function itself was used, making the construction less canonical.

\subsubsection*{The choice of $N$}
By the results of Kouchnirenko \cite{Kouchnirenko76}, the subspace $\Afr$ of $\CC[u,u^{-1}]$ generated by the monomials $u^{\alpha_i}$, where $(\alpha_i)_{i=1,\dots,r}$ denote the integral points of $\ZZ^n$ contained in the interior of the Newton polyhedron of $f$, injects to the Jacobi space $\CC[u,u^{-1}]/(\partial f)$.

We consider the deformation
\[
F(u,x)=f(u)+\sum_{i=1}^rx_iu^{\alpha_i}
\]
parametrized by the affine space $\Afr$. Clearly, for any $x\in\Afr$, the function $F_x(u)\defin F(u,x)$ is convenient and nondegenerate. If $f=F_0$ satisfies \eqref{cond:semissimple}, there exists a Zariski dense open subset $N$ of $\Afr$ such that $F_x$ satisfies \eqref{cond:semissimple} for any $x\in N$. If $F_0$ satisfies \eqref{cond:strongconv}, then so does $F_x$ for any $x\in N\defin\Afr$.

\subsubsection*{The Saito structure on $N$}
We will first describe the canonical Saito structure parametrized by the complex manifold $\Afr$ attached to $F$.

We denote by $\Omega^k(U)[x]$ the space of algebraic differential $k$-forms on $U\times\Afr$ relative to $\Afr$ (\ie with no $dx$) and by $d_u$ the relative differential. The Brieskorn lattice $G_{F,0}$ of $F$ is the quotient
\[
\Omega^n(U)[x,\hb]\big/(\hb d_u-d_uF\wedge)\Omega^{n-1}(U)[x,\hb].
\]
Under our assumption on $f$, this is a free $\CC[\hb,x_1,\dots,x_n]$-module of rank $\mu(f)$, equipped with a connection $\nabla$ having Poincar\'e rank one along $\hb=0$ (\cf \cite{Douai05}). The connection is determined by its action on the class $[\eta]$ of differential forms $\eta\in\Omega^n(U)$:
\[
\hb^2\nabla_{\partial_\hb}[\eta]=[F\eta],\quad\hb\nabla_{\partial_{x_i}}[\eta]=-[F'_{x_i}\eta]\quad(i=1,\dots,r).
\]

For any $x^o\in\Afr$ (for instance $x^o=0$), let us consider the corresponding object attached to the function $F^o(u)=F(u,x^o)$. By using a construction of M\ptbl Saito together with Hodge theory for the function $F^o$ (\cf \cite{MSaito89,Bibi96bb,Hertling00,D-S02a}) one can define in a canonical way a trivialization $G_{F^o,0}\simeq\CC[\hb]^{\mu(f)}$ so that the connection $\nabla^o$ takes the form $d+(\cU^o/\hb-\cV^o)d\hb/\hb$, $\cV^o$ is semi-simple and its eigenvalues form the \emph{spectrum at infinity} of $F^o$ translated by $-n/2$ (that is, centered at~$0$), together with a canonical pairing $\cP^o$ obtained by microlocalization from the Poincar\'e duality in the fibres of $F^o$.

A\ptbl Douai moreover shows (\cf \cite{Douai05} after the proof of Prop\ptbl3.1.2 and \cite[Th\ptbl3.3.1]{Douai07}) that this construction can be made for $F$ over $\Afr$, is compatible with base change, and gives a matrix for the connection $\nabla$ as in \eqref{eq:nablafnabla} and a pairing $\cP$ as in \eqref{eq:cP}.

\subsubsection*{The pre-primitive section}
Now, the (constant w.r.t.~$x$) complex volume form $\omega=\dfrac{du_1}{u_1}\wedge\cdots\wedge\dfrac{du_n}{u_n}$ defines all along $\Afr$ a real homogeneous pre-primitive form (it corresponds to the eigenvalue $-n/2$ of $\cV$).

The choice of $N$ (possibly distinct from $\Afr$) is made to apply, in Case \eqref{cond:semissimple}, the theorem of Malgrange at each point of $N$. In Case \eqref{cond:strongconv}, we can apply the result of \cite{H-M04} all along $\Afr$.

\subsection{Twistor structures associated to tame functions}\label{subsec:tame}
We now consider more generally a $n$-dimensional smooth affine variety $U$ and a regular function $f:U\to\CC$, that we assume to be ``cohomologically tame'': this property means that the critical points of $f$ in $U$ are isolated and that no critical value of $f$ comes from a critical point at infinity (see \eg \cite{Bibi96bb}). For instance, if $U$ is a torus and $f$ is a Laurent polynomial, tameness is implied by the property of being convenient and nondegenerate. We will attach below to $f$ a canonical twistor structure $(\cH^{\prime o},\nabla^o,\cCS^o)$ of weight~$0$.

Let us first define $(\cH^{\prime o},\nabla^o)$. The bundle $\cH^{\prime o}$ is the analytization of the Brieskorn lattice
\[
\Omega^n(U)[\hb]/(\hb d-df\wedge)\Omega^{n-1}(U)[\hb],
\]
which is known to be a free $\CC[\hb]$-module of finite rank $\mu(f)$ (\cf \eg \cite{Bibi96bb}) equipped with the connection $\nabla^o$ defined by
\[
\hb^2\nabla^o_{\partial_\hb}\big[\textstyle\sum_k\omega_k\hb^k\big]= \big[\sum_kk\omega_k\hb^{k+1}+\sum_kf\omega_k\hb^k\big].
\]

On the other hand, it is classical (after the work of F\ptbl Pham \cite{Pham85b}) that the local system $\cL^o$ is identified to the locally constant sheaf of Lefschetz co-thimbles
\[
H^n_{\Phi_\hb}(U,\CC),
\]
where $H^*_{\Phi_\hb}$ denotes the cohomology with support in the family $\Phi_\hb$ of closed sets in~$U$ on which $\reel(f(u_1,\dots,u_n)/\hb)\leq c<0$.

There is a natural intersection pairing (Poincar\'e duality pairing made sesquilinear)
\begin{equation}\label{eq:Poincare}
\wh P_\hb:H^n_{\Phi_\hb}(U,\CC)\otimes \ov{H^n_{\Phi_{-\hb}}(U,\CC)}\to \CC.
\end{equation}

\begin{definition}[C\ptbl Hertling]\label{def:cCS}
The twistor structure canonically attached to $f$ is $(\cH^{\prime o},\nabla^o,\cCS^o)$, where $(\cH^{\prime o},\nabla^o)$ is the analytization of the Brieskorn lattice of~$f$ and $\cCS^o\defin\dfrac{(-1)^{n(n-1)/2}}{(2\pi i)^n}\wh P$.
\end{definition}

\begin{theoreme}[{\cf \cite[Th\ptbl4.10]{Bibi05}}]\label{th:positivityftame}
The twistor structure $(\cH^{\prime o},\nabla^o,\cCS^o)$ attached to~$f$ is pure of weight~$0$ and polarized.
\end{theoreme}

\begin{proof}[Indication for the proof]\mbox{}

\subsubsection*{Reduction of the problem to dimension one}
How can one prove such a positivity statement? One should start with a variation of polarized twistor structure of weight~$0$ and get our twistor structure by a natural operation from the previous one. Here is an example of such a result:

\begin{HStheoreme}
Given a variation of polarized twistor structure of weight~$0$ on a compact K\"ahler manifold $X$, its de~Rham cohomology carries a polarized twistor structure (of~some weight).
\end{HStheoreme}

The main ingredient in the proof of the previous theorem is the fact that, for any $\hb\in\CC$, the Laplace operator $\Delta_\hb$ relative to the operator $D_\hb=d''+\Phi+\hb(D'+\Phi^\dag)$ and the K\"ahler metric is essentially constant: $\Delta_\hb=(1+|\hb|^2)\Delta_0$ (this is the analogue of the classical K\"ahler identity $\Delta_d=2\Delta_{d'}=2\Delta_{d''}$). Hence, the space of harmonic sections does not depend on $\hb$. This will give the pure weight~$0$ property. The positivity is obtained by a standard argument of Hodge theory, on primitive sections first.

In Theorem \ref{th:positivityftame}, one can use a similar argument (taking direct image by the constant map, that is, taking de~Rham cohomology): one can obtain $(\cH^{\prime o},\nabla^o)$ by
\begin{enumerate}
\item
considering the trivial variation of twistor structure $(\cO_U[\hb],\hb d)$ (the Higgs field is equal to $0$),
\item
twisting it by $e^{-f/\hb}$, that is, adding a new Higgs field $\Phi=-df$,
\item
and taking the de~Rham cohomology of this new variation.
\end{enumerate}
The corresponding operator $D_\hb$ is $d''-df+\hb(d'-d\ov f)$. We are now faced with two problems: $U$~is noncompact and $f$ is not bounded on $U$ (so that $e^{-f}$ can have an exponential growth). The Hodge theory for the corresponding Laplacian can be difficult to develop (although it has been developed in some special cases).

Instead, we use \emph{Horatio's method}: if we face numerous enemies, we fake escaping by running fast, then kill the enemy running faster when he reaches us, then kill the next one, etc. Here, we escape by falling down along the fibres of $f:U\to\Afu$.

Let $t$ be the coordinate on $\Afu$. The Gauss-Manin connection of $f$ gives bundles with connection on $\Afu\moins\{\text{critical values of }f\}$. The interesting bundle $V^\an$ has fibre $H^{n-1}(f^{-1}(t),\CC)$. It underlies a variation of mixed Hodge structure (M\ptbl Saito). The assumption made on $f$ (cohomological tameness) implies that this mixed Hodge structure is an extension of pure Hodge structures for which one subquotient is a variation of polarized Hodge structure (whose generic fibre is identified to the intersection cohomology of a suitable compactification of $f^{-1}(t)$) and any other quotient is a trivial variation of Hodge structure on $\Afu$.

The variation of polarized Hodge structure induces a variation of polarized twistor structure of the same weight.

\subsubsection*{End of the proof of the theorem}
We will give the main steps (\cf \cite{Bibi05} for details).

Starting from $f:U\to\Afu$, we consider the Gauss-Manin system $M$ of $f$: $M=\Omega^n(U)[\partial_t]/(d-\partial_tdf\wedge)\Omega^{n-1}(U)[\partial_t]$ with its natural structure of $\CC[t]\langle\partial_t\rangle$-module. According to a general result of M\ptbl Saito \cite{MSaito87}, it underlies a~mixed Hodge module on $\Afu$ (hence it is equipped with a Hodge filtration $F_\bbullet^\rH M$). This mixed Hodge module has a subquotient $M_{!*}$ (in the category of mixed Hodge modules) which is a polarized pure Hodge module of weight $n$, and $M$ is an extension of $M_{!*}$ by \emph{constant} Hodge modules (\ie isomorphic to powers of $\CC[t]$ with its natural $\CC[t]\langle\partial_t\rangle$-action and trivial Hodge filtration up to a shift): this is a consequence of the tameness assumption.

Let us set $F_\ell M_{!*}\defin F_{\ell-n}^\rH M_{!*}$ and let us choose an index $\ell_0$ such that $M_{!*}/F_{\ell_0}M_{!*}$ is supported on a finite set of points (this is always possible by definition of a good filtration). We use $F_{\ell_0}M_{!*}$ as a submodule $L$ considered in \S\ref{subsec:Laplace}.

We note that, considering the localized Fourier-Laplace transforms of $M$ and $M_{!*}$, we have $G=G_{!*}$, as the localized Fourier-Laplace transform of $\CC[t]$ is zero. But in a less trivial way, if we denote by $G_{!*,0}$ the $\CC[\hb]$-module $\hb^{\ell_0}G_{!*,0}^{(L)}$ (as defined in \S\ref{subsec:Laplace}) and keep the notation $G_0$ for the Brieskorn lattice of $f$, we have (\cf \cite[Lemma 4.7]{Bibi05}) $G_{!*,0}=G_0$. In other words, this equality gives a definition of the Brieskorn lattice purely in terms of Hodge theory.

We now apply the construction of \S\ref{sec:Fourier-Laplace} to $M_{!*}$ and $G_{!*,0}$. Theorem \ref{th:FourierHodge} shows that it is a polarized pure twistor structure of weight~$0$.

It remains to identify the $\cC_{!*,\bS}$ with $\cCS^o$ given in Definition \ref{def:cCS}. This is done by analyzing the definition of the polarization given by M\ptbl Saito \cite{MSaito86}.
\end{proof}

\subsection{The canonical harmonic Frobenius manifold structure}
We are now in position to state the main theorem. Let $f$ be a convenient and nondegenerate Laurent polynomial. Assume that $f$ satisfies one of the conditions \eqref{cond:semissimple} or \eqref{cond:strongconv} and let $M$ be the germ along $N$ of the $\CC$-vector space $\CC[u,u^{-1}]/(\partial f)$, equipped with its canonical Frobenius manifold structure, as explained in \S\ref{subsec:Douai}. For any $x^o\in N$, the twistor structure given by Theorem \ref{th:positivityftame} at $x^o$ extends locally, by using the flat extension of the pairing $\cCS^o$ of Definition \ref{def:cCS} and, according to Corollary \ref{cor:deformation}, it defines a positive definite Hermitian form on the Saito structure attached to $F$ in some neighbourhood of $x^o$. Transporting this Hermitian form by the infinitesimal period mapping associated to the primitive real homogeneous section $\dfrac{du_1}{u_1}\wedge\cdots\wedge\dfrac{du_n}{u_n}$ endows the germ $(M,x^o)$ with a canonical harmonic Frobenius manifold structure (\cf Corollary \ref{cor:harmonic}) with a positive definite harmonic metric.

\begin{theoreme}\label{th:main}
Under these conditions, the canonical harmonic Frobenius manifold structures on $(M,x^o)$ glue together when $x^o$ varies in $N$ and define a canonical harmonic Frobenius manifold structure with a positive definite harmonic metric on the germ $(M,N)$.
\end{theoreme}

\begin{proof}
The point is to show that the flat extension $\cC_\bS$ of $\cC_\bS^o$ restricts, up to the factor $(-1)^{n(n-1)/2}/(2\pi i)^n$, to the pairing \eqref{eq:Poincare} defined from $F_x$ by Poincaré duality. In other words, one should give a definition of the pairing \eqref{eq:Poincare} with parameter and show the compatibility with base change. Let us set $\wt F(u,x)=(F(u,x),x):U\times\Afr\to\Afu\times\Afr$. Then $R\wt F_!\CC_{U\times\Afr}$ is a constructible complex on $\Afu\times\Afr$ which is compatible with base change with respect to $x\in\Afr$. It is also Poincaré-Verdier dual to $R\wt F_*\CC_{U\times\Afr}$ up to a shift.

\begin{assertion*}
Given any $x^o\in\Afr$, there exists a disc $D^o$ in $\Afu$ such that the cohomology sheaves of $R\wt F_!\CC_{U\times\Afr}$ (or $R\wt F_*\CC_{U\times\Afr}$) are local systems on $(\Afu\moins D^o)\times\nb(x^o)$.
\end{assertion*}

Let us postpone the proof of the assertion. The partial Laplace transform of $R\wt F_!\CC_{U\times\Afr}$ (or $R\wt F_*\CC_{U\times\Afr}$) with respect to $\Afu$ gives then, by the same computations as in \cite[\S1b]{Bibi05}, a constructible complex on $(\Afu_\tau\moins\{0\})\times\Afr$ whose cohomologies are local systems compatible with base change with respect to $x\in\Afr$, hence there is only one nonzero cohomology local system.

Note that the pairing \eqref{eq:Poincare} at any $x\in\Afr$ comes by Laplace transform from the Poincaré duality paring between $RF_{x,!}\CC_U$ and $RF_{x,*}\CC_U$ made sesquilinear, and also from its restriction between $RF_{x,!}\CC_U$ and $RF_{x,!}\CC_U$ because $RF_{x,*}\CC_U$ and $RF_{x,!}\CC_U$ have the same Laplace transform. It is therefore obtained by base change from the partial Laplace transform of the restriction of the Poincaré duality pairing (made sesquilinear) between $R\wt F_!\CC_{U\times\Afr}$ and itself. By definition, the latter object is the flat extension of $\wh P_{x}$.
\end{proof}

\begin{proof}[Sketch of a proof for the assertion]
Assume first that there exists a basis of $\ZZ^n$ contained in the interior of the Newton polygon of $f$. We will work with the corresponding coordinates $u_i$ and argue as in \cite{Broughton88}. The Milnor number of any $F_x$ is constant with respect to $x\in\Afr$ and so is the Milnor number of $F_x^w\defin F_x+\sum_iw_iu_i$ (note that this deformation of $f$, taken with parameters $x$ and $w$, is somewhat redundant). As a consequence, given a compact neighbourhood $\ov B$ of the critical points of $f$, there exists a compact neighbourhood $\nb(x^o)$ of $x^o$ such that the pull-back by $(\partial F/\partial u_1,\dots,\partial F/\partial u_n):(\CC^*)^n\times\nb(x^o)\to\CC^n$ of some compact neighbourhood of the origin is contained in $\ov B\times\nb(x^o)$. In other words, $|\partial F/\partial u_i|$ are uniformly bounded from below on $[(\CC^*)^n\moins B]\times\nb(x^o)$. Moreover, as $\ov B$ is compact, the previous set is contained in $\wt F^{-1}\big((\Afu\moins D^o)\times\nb(x^o)\big)$ for a suitable disc $D^o\subset\Afu$. Now the argument used in \cite[p\ptbl230]{Broughton88} to show the fibration property extends to the present situation and shows that $\wt F$ is a locally trivial topological fibration above $(\Afu\moins D^o)\times\nb(x^o)$, hence the result.

If such a basis does not exist, we first perform a covering $\rho:V\to U$ of the kind $u_i=v_i^m$, for $m$ large enough to reduce to the previous case. The assertion then holds for $R\wt F_!\rho_*\CC_{V\times\Afr}$ and we take invariants with respect to the covering group action to obtain the result.
\end{proof}

\providecommand{\bysame}{\leavevmode\hbox to3em{\hrulefill}\thinspace}
\providecommand{\MR}{\relax\ifhmode\unskip\space\fi MR }
\providecommand{\MRhref}[2]{%
  \href{http://www.ams.org/mathscinet-getitem?mr=#1}{#2}
}
\providecommand{\href}[2]{#2}


\begin{thebibliography}{10}

\bibitem{Broughton88}
{\relax S.A}.~Broughton, \emph{Milnor number and the topology of polynomial
  hypersurfaces}, Inventiones Math. \textbf{92} (1988), 217--241.

\bibitem{Corlette88}
K.~Corlette, \emph{Flat {$G$}-bundles with canonical metrics}, J.~Differential
  Geom. \textbf{28} (1988), 361--382.

\bibitem{Douai05}
A.~Douai, \emph{Construction de vari\'et\'es de {F}robenius via les polyn\^omes
  de {L}aurent: une autre approche}, Singularit\'es, Rev. Inst. \'Elie Cartan,
  vol.~18, Univ. Nancy, Nancy, 2005, pp.~105--123.

\bibitem{Douai07}
\bysame, \emph{{A canonical Frobenius structure}}, arXiv: \url{0709.0186v1},
  2007.

\bibitem{D-S02a}
A.~Douai and C.~Sabbah, \emph{{Gauss-Manin systems, Brieskorn lattices and
  Frobenius structures (I)}}, Ann. Inst. Fourier (Grenoble) \textbf{53} (2003),
  no.~4, 1055--1116.

\bibitem{Hertling00}
C.~Hertling, \emph{{Frobenius manifolds and moduli spaces for singularities}},
  Cambridge Tracts in Mathematics, vol. 151, Cambridge University Press, 2002.

\bibitem{Hertling01}
\bysame, \emph{{$tt^*$ geometry, Frobenius manifolds, their connections, and
  the construction for singularities}}, J.~reine angew. Math. \textbf{555}
  (2003), 77--161.

\bibitem{H-M04}
C.~Hertling and {\relax Yu.I}.~Manin, \emph{{Unfoldings of meromorphic
  connections and a construction of Frobenius manifolds}}, {Frobenius manifolds
  (Quantum cohomology and singularities)} (C.~Hertling and M.~Marcolli, eds.),
  Aspects of Mathematics, vol. E~36, Vieweg, 2004, pp.~113--144.

\bibitem{H-S06}
C.~Hertling and {\relax Ch}.~Sevenheck, \emph{{Nilpotent orbits of a
  generalization of Hodge structures}}, J.~reine angew. Math. \textbf{609}
  (2007), 23--80, arXiv: \url{math.AG/0603564}.

\bibitem{Kashiwara86}
M.~Kashiwara, \emph{Regular holonomic {$\mathcal D$}-modules and distributions
  on complex manifolds}, Complex analytic singularities, Adv. Stud. Pure Math.,
  vol.~8, North-Holland, Amsterdam, 1987, pp.~199--206.

\bibitem{Kouchnirenko76}
{\relax A.G}.~Kouchnirenko, \emph{{Poly\`edres de Newton et nombres de
  Milnor}}, Invent. Math. \textbf{32} (1976), 1--31.

\bibitem{Leiterer90}
J.~Leiterer, \emph{{Holomorphic vector bundles and the Oka-Grauert principle}},
  {Several complex variables. IV. Algebraic aspects of complex analysis},
  Encycl. Math. Sci., vol.~10, Springer-Verlag, 1990, pp.~63--103.

\bibitem{Malgrange83db}
B.~Malgrange, \emph{{Sur les d\'eformations isomonodromiques, I, II}},
  {S\'eminaire E.N.S. Math\'ematique et Physique} (L.~Boutet~{de Monvel},
  A.~Douady, and J.-L. Verdier, eds.), Progress in Math., vol.~37,
  Birkh{\"a}user, Basel, Boston, 1983, pp.~401--438.

\bibitem{Malgrange91}
\bysame, \emph{{\'E}quations diff\'erentielles {\`a} coefficients polynomiaux},
  Progress in Math., vol.~96, Birkh{\"a}user, Basel, Boston, 1991.

\bibitem{Mochizuki07}
T.~Mochizuki, \emph{{Asymptotic behaviour of tame harmonic bundles and an
  application to pure twistor $D$-modules}}, vol. 185, Mem. Amer. Math. Soc.,
  no. 869-870, American Mathematical Society, Providence, RI, 2007.

\bibitem{Pham85b}
F.~Pham, \emph{{La descente des cols par les onglets de Lefschetz avec vues sur
  Gauss-Manin}}, {Syst\`emes diff\'erentiels et singularit\'es (Luminy, 1983)}
  (A.~Galligo, J.-M. Granger, and {\relax Ph}.~Maisonobe, eds.),
  Ast{\'e}risque, vol. 130, Soci{\'e}t{\'e} Math{\'e}matique de France, 1985,
  pp.~11--47.

\bibitem{Bibi00}
C.~Sabbah, \emph{{D\'eformations isomonodromiques et vari\'et\'es de
  Frobenius}}, Savoirs Actuels, CNRS~{\'E}ditions \& EDP~Sciences, Paris, 2002,
  English Transl.: Universitext, Springer \& EDP~Sciences, 2007.

\bibitem{Bibi99}
\bysame, \emph{{Vanishing cycles and Hermitian duality}}, Proc. Steklov Inst.
  Math. \textbf{238} (2002), 194--214.

\bibitem{Bibi04}
\bysame, \emph{{Fourier-Laplace transform of irreducible regular differential
  systems on the Riemann sphere}}, Russian Math. Surveys \textbf{59} (2004),
  no.~6, 1165--1180, erratum:
  \url{http://math.polytechnique.fr/~sabbah/sabbah_Fourier-irred_err.pdf}.

\bibitem{Bibi01c}
\bysame, \emph{{Polarizable twistor $\mathcal{D}$-modules}}, Ast{\'e}risque,
  vol. 300, Soci{\'e}t{\'e} Math{\'e}matique de France, Paris, 2005.

\bibitem{Bibi96bb}
\bysame, \emph{Hypergeometric periods for a tame polynomial}, Portugal. Math.
  \textbf{63} (2006), no.~2, 173--226, arXiv: \url{math.AG/9805077}.

\bibitem{Bibi05}
\bysame, \emph{{Fourier-Laplace transform of a variation of polarized complex
  Hodge structure}}, J.~reine angew. Math. (to appear), arXiv:
  \url{math.AG/0508551}, 36 pages.

\bibitem{KSaito83b}
K.~Saito, \emph{Period mapping associated to a primitive form}, Publ. RIMS,
  Kyoto Univ. \textbf{19} (1983), 1231--1264.

\bibitem{MSaito86}
M.~Saito, \emph{Modules de {Hodge} polarisables}, Publ. RIMS, Kyoto Univ.
  \textbf{24} (1988), 849--995.

\bibitem{MSaito89}
\bysame, \emph{{On the structure of Brieskorn lattices}}, Ann. Inst. Fourier
  (Grenoble) \textbf{39} (1989), 27--72.

\bibitem{MSaito87}
\bysame, \emph{{Mixed {Hodge} Modules}}, Publ. RIMS, Kyoto Univ. \textbf{26}
  (1990), 221--333.

\bibitem{Schmid73}
W.~Schmid, \emph{Variation of {Hodge} structure: the singularities of the
  period mapping}, Invent. Math. \textbf{22} (1973), 211--319.

\bibitem{Simpson88}
C.~Simpson, \emph{Constructing variations of {Hodge} structure using
  {Yang-Mills} theory and applications to uniformization}, J.~Amer. Math. Soc.
  \textbf{1} (1988), 867--918.

\bibitem{Simpson90}
\bysame, \emph{Harmonic bundles on noncompact curves}, J.~Amer. Math. Soc.
  \textbf{3} (1990), 713--770.

\bibitem{Simpson92}
\bysame, \emph{Higgs bundles and local systems}, Publ. Math. Inst. Hautes
  {\'E}tudes Sci. \textbf{75} (1992), 5--95.

\bibitem{Simpson97}
\bysame, \emph{Mixed twistor structures}, Pr\'epublication Universit\'e de
  Toulouse \& arXiv: \url{math.AG/9705006}, 1997.

\bibitem{Szabo04}
S.~Szabo, \emph{{Nahm transform of meromorphic integrable connections on the
  Riemann sphere}}, Ph.D. thesis, Universit\'e Louis Pasteur, Strasbourg,
  juillet 2005, arXiv: \url{0511471v1}.

\bibitem{Takahashi04}
A.~Takahashi, \emph{{tt$^*$ geometry of rank two}}, Internat. Math. Res.
  Notices (2004), no.~22, 1099--1114.

\bibitem{Zucker79}
S.~Zucker, \emph{{Hodge theory with degenerating coefficients:
  {$L_2$}-cohomology in the Poincar\'e metric}}, Ann. of Math. \textbf{109}
  (1979), 415--476.

\end{thebibliography}
\end{document}